\documentclass[10pt]{article}

\usepackage{amsmath}
\usepackage{amsthm}
\usepackage{amssymb}
\usepackage{amscd}
\usepackage{amsfonts}  
\usepackage[all,dvips]{xy}
\usepackage{makeidx}     
\usepackage{graphicx}    
\usepackage{multicol}    

\DeclareMathAlphabet\EuR{U}{eur}{m}{n}
\SetMathAlphabet\EuR{bold}{U}{eur}{b}{n}

\begin{document}


\newtheorem{theorem}{Theorem}[section]
\newtheorem{lemma}[theorem]{Lemma}
\newtheorem{proposition}[theorem]{Proposition}
\newtheorem{definition}[theorem]{Definition}
\newtheorem{example}[theorem]{Example}
\newtheorem{remark}[theorem]{Remark}
\newtheorem{corollary}[theorem]{Corollary}
\newtheorem{conjecture}[theorem]{Conjecture}
\newtheorem{convention}[theorem]{Convention}
\newtheorem{problem}[theorem]{Problem}
{\catcode`@=11\global\let\c@equation=\c@theorem}
\renewcommand{\theequation}{\thetheorem}

\renewcommand{\theenumi}{\alph{enumi}}
\renewcommand{\theenumii}{\roman{enumii}}
\renewcommand{\labelenumi}{(\theenumi)}
\renewcommand{\labelenumii}{(\theenumii)}

\makeatletter
\renewcommand{\@seccntformat}[1]{\csname the#1\endcsname.\hspace{1em}}
\makeatother
\newcommand{\tit}[2]{\begin{bf} \begin{center} \begin{Large}
\section{#1}
\label{sec:#2}
\end{Large}\end{center}\end{bf}
\nopagebreak}


\newcommand{\squarematrix}[4]{\left( \begin{array}{cc} #1 & #2 \\ #3 &
#4
\end{array} \right)}
\newcommand{\smallmx}[4]{\mbox{\begin{scriptsize}$\squarematrix{#1}{#2}
        {#3}{#4}$\end{scriptsize}}}

\newcommand{\comsquare}[8]{
\begin{center}
$\begin{CD}
#1 @>#2>> #3\\
@V{#4}VV @VV{#5}V\\
#6 @>>#7> #8
\end{CD}$
\end{center}}

\newcommand{\comsquarewithout}[8]{
\begin{CD}
#1 @>#2>> #3\\
@V{#4}VV @VV{#5}V\\
#6 @>>#7> #8
\end{CD}}


\let\sect=\S
\newcommand{\curs}{\EuR}
\newcommand{\CHAINCOMPLEXES}{\curs{CHCOM}}
\newcommand{\GROUPOIDS}{\curs{GROUPOIDS}}
\newcommand{\PAIRS}{\curs{PAIRS}}
\newcommand{\FGINJ}{\curs{FGINJ}}
\newcommand{\Or}{\curs{Or}}
\newcommand{\SPACES}{\curs{SPACES}}
\newcommand{\SPECTRA}{\curs{SPECTRA}}
\newcommand{\Sub}{\curs{Sub}}


\newcommand{\bbC}{{\mathbb C}}
\newcommand{\bbH}{{\mathbb H}}
\newcommand{\bbI}{{\mathbb I}}
\newcommand{\bbK}{{\mathbb K}}
\newcommand{\bbKO}{\mathbb{KO}}
\newcommand{\bbN}{{\mathbb N}}
\newcommand{\bbP}{{\mathbb P}}
\newcommand{\bbQ}{{\mathbb Q}}
\newcommand{\bbR}{{\mathbb R}}
\newcommand{\bbZ}{{\mathbb Z}}

\newcommand{\calbh}{{\mathcal B}{\mathcal H}}
\newcommand{\calc}{{\mathcal C}}
\newcommand{\cald}{{\mathcal D}}
\newcommand{\cale}{{\mathcal E}}
\newcommand{\calf}{{\mathcal F}}
\newcommand{\calg}{{\mathcal G}}
\newcommand{\calh}{{\mathcal H}}
\newcommand{\calk}{{\mathcal K}}
\newcommand{\call}{{\mathcal L}}
\newcommand{\caln}{{\mathcal N}}
\newcommand{\cals}{{\mathcal S}}

\newcommand{\bfA}{\ensuremath{\mathbf{A}}}
\newcommand{\bfa}{\ensuremath{\mathbf{a}}}
\newcommand{\bfb}{\ensuremath{\mathbf{b}}}
\newcommand{\bfCTR}{\ensuremath{\mathbf{CTR}}}
\newcommand{\bfDtr}{\ensuremath{\mathbf{Dtr}}}
\newcommand{\bfE}{\ensuremath{\mathbf{E}}}
\newcommand{\bff}{\ensuremath{\mathbf{f}}}
\newcommand{\bfF}{\ensuremath{\mathbf{F}}}
\newcommand{\bfg}{\ensuremath{\mathbf{g}}}
\newcommand{\bfHH}{\ensuremath{\mathbf{HH}}}
\newcommand{\bfHPB}{\ensuremath{\mathbf{HPB}}}
\newcommand{\bfI}{\ensuremath{\mathbf{I}}}
\newcommand{\bfi}{\ensuremath{\mathbf{i}}}
\newcommand{\bfK}{\ensuremath{\mathbf{K}}}
\newcommand{\bfL}{\ensuremath{\mathbf{L}}}
\newcommand{\bfr}{\ensuremath{\mathbf{r}}}
\newcommand{\bfs}{\ensuremath{\mathbf{s}}}
\newcommand{\bfS}{\ensuremath{\mathbf{S}}}
\newcommand{\bft}{\ensuremath{\mathbf{t}}}
\newcommand{\bfT}{\ensuremath{\mathbf{T}}}
\newcommand{\bfTC}{\ensuremath{\mathbf{TC}}}
\newcommand{\bfU}{\ensuremath{\mathbf{U}}}
\newcommand{\bfu}{\ensuremath{\mathbf{u}}}
\newcommand{\bfv}{\ensuremath{\mathbf{v}}}
\newcommand{\bfw}{\ensuremath{\mathbf{w}}}


\newcommand{\Arf}{\operatorname{Arf}}
\newcommand{\aut}{\operatorname{aut}}
\newcommand{\Bor}{\operatorname{Bor}}
\newcommand{\ch}{\operatorname{ch}}
\newcommand{\class}{\operatorname{class}}
\newcommand{\cok}{\operatorname{coker}}
\newcommand{\cone}{\operatorname{cone}}
\newcommand{\colim}{\operatorname{colim}}
\newcommand{\con}{\operatorname{con}}
\newcommand{\conhom}{\operatorname{conhom}}
\newcommand{\cyclic}{\operatorname{cyclic}}
\newcommand{\ev}{\operatorname{ev}}
\newcommand{\ext}{\operatorname{ext}}
\newcommand{\fr}{\operatorname{fr}}
\newcommand{\Gen}{\operatorname{Gen}}
\newcommand{\hoaut}{\operatorname{ho-aut}}
\newcommand{\hur}{\operatorname{hur}}
\newcommand{\im}{\operatorname{im}}
\newcommand{\inj}{\operatorname{inj}}
\newcommand{\id}{\operatorname{id}}
\newcommand{\infl}{\operatorname{Infl}}
\newcommand{\ind}{\operatorname{ind}}
\newcommand{\Inn}{\operatorname{Inn}}
\newcommand{\Irr}{\operatorname{Irr}}
\newcommand{\Is}{\operatorname{Is}}
\newcommand{\ks}{\operatorname{ks}}
\newcommand{\map}{\operatorname{map}}
\newcommand{\MOD}{\operatorname{MOD}}
\newcommand{\mor}{\operatorname{mor}}
\newcommand{\Ob}{\operatorname{Ob}}
\newcommand{\op}{\operatorname{op}}
\newcommand{\pr}{\operatorname{pr}}
\newcommand{\point}{\operatorname{pt.}}
\newcommand{\PT}{\operatorname{PT}}
\newcommand{\Rat}{\operatorname{Rat}}
\newcommand{\res}{\operatorname{res}}
\newcommand{\sign}{\operatorname{sign}}
\newcommand{\Th}{\operatorname{Th}}
\newcommand{\topo}{\operatorname{top}}
\newcommand{\tors}{\operatorname{tors}}
\newcommand{\Wh}{\operatorname{Wh}}

\newcommand{\pt}{\{\bullet\}}

\newcounter{commentcounter}
\newcommand{\commentw}[1]                      
{\stepcounter{commentcounter}
{\bf Comment \arabic{commentcounter} (by W.)}: {\ttfamily #1} }

\newcommand{\commentm}[1]                      
{\stepcounter{commentcounter}
{\bf Comment \arabic{commentcounter} by M.}: {\ttfamily #1} }

\hyphenation{equi-va-riant}
\hyphenation{mani-fold}
\hyphenation{di-men-sion}
\hyphenation{di-men-sio-nal}
\hyphenation{Con-jec-ture}

\title{Topological rigidity for non-aspherical manifolds\\
by \\
M. Kreck and W. L\"uck}
\maketitle



\typeout{-----------------------  Abstract  ------------------------}

\begin{abstract}
The Borel Conjecture predicts that closed aspherical manifolds are topological rigid.
We want to investigate when a non-aspherical oriented connected closed manifold $M$ is topological rigid in the following sense.
If $f \colon N \to M$ is an orientation preserving homotopy equivalence with a closed oriented manifold as target,
then there is an orientation preserving homeomorphism $h \colon N \to M$ such that $h$ and $f$ induce up to conjugation
the same maps on the fundamental groups. We call such manifolds 
\emph{Borel manifolds}. 
We give partial answers to this questions for $S^k \times S^d$,
for sphere bundles over aspherical closed manifolds of dimension $\le 3$ and for $3$-manifolds with torsionfree fundamental groups.
We show that this rigidity is inherited under connected sums in dimensions $\ge 5$.
We also classify manifolds of dimension $5$ or $6$ whose fundamental group is the one of a surface and whose second
homotopy group is trivial.\\[1mm]
\end{abstract}


\typeout{--------------------------------   Section 0: Introduction and Statement of Results ------------------------------------}

\setcounter{section}{-1}

\section{Introduction and Statement of Results}
\label{sec:Introduction_and_Statement_of_Results}

In this paper we study the question which non-aspherical oriented closed connected topological manifolds
are topological rigid. Recall that the Borel Conjecture predicts that every aspherical closed
topological manifold is topological rigid in the sense
that every homotopy equivalence of such manifolds is homotopic to a homeomorphism.
We focus on the following two problems which we will describe next.

We say that two maps $f,g \colon X \to Y$ of path-connected spaces
\emph{induce the same map on the fundamental groups up to conjugation}
if for one (and hence all base points) $x \in X$ there exists a path $w$ from $f(x)$ to $g(x)$ such that
for the group isomorphism $t_w \colon \pi_1(Y,f(x)) \to \pi_1(Y,g(y))$
which sends the class of a loop $v$ to the class of the loop
$w^- \ast v \ast w$ we get $\pi_1(g,x) = t_w \circ \pi_1(f,x)$.
Homotopic maps induce the same map on the fundamental groups up to conjugation.

\begin{convention}
\label{con:manifold}
Manifold will always mean connected oriented closed topological manifold unless stated explicitly differently.
\end{convention}

\begin{definition}[Borel-manifold]\label{def: Borel manifold} \em
A manifold $M$ is called a \emph{Borel manifold}
if for any orientation preserving homotopy equivalence $f \colon N \to M$ of manifolds
there exists an orientation preserving homeomorphism $h \colon N \to M$ such that
$f$ and $h$ induce the same map on the fundamental groups up to
conjugation.
It is called a \emph{strong Borel manifold}
if every orientation preserving homotopy equivalence $f \colon N \to M$ of manifolds
is homotopic to a homeomorphism $h \colon N \to M$.
\em
\end{definition}

\begin{remark}[Relation to the Borel Conjecture] \label{rem:Relation to the Borel Conjecture}
\em
If $M$ is aspherical, two homotopy equivalences $f,g \colon N \to M$ are homotopic if and only if they induce
the same map on the fundamental groups up to conjugation. Hence an aspherical manifold
$M$ is a Borel manifold if and only if
every homotopy equivalence $f \colon N \to M$ of manifolds is
homotopic to a homeomorphism. This is the precise statement of the
\emph{Borel Conjecture} for aspherical manifolds.
Hence the Borel Conjecture can be rephrased as the statement that every aspherical manifold
$M$ is a Borel manifold, or equivalently, is a strong Borel manifold.
More information on the Borel Conjecture can be found
for instance in \cite{Farrell(1996)},
\cite{Farrell(2002)}, \cite{Farrell-Jones(1989)}, \cite{Farrell-Jones(1990)}, \cite {Farrell-Jones(1993c)},
\cite{Farrell-Jones(1998)},
\cite{Ferry-Ranicki-Rosenberg(1995)}, \cite{Kreck-Lueck(2005)}, \cite{Lueck(2002c)},
\cite{Lueck-Reich(2005)}.
\end{remark}

\begin{remark}[Relation to the Poincar\'e Conjecture] \label{rem:Relation to the Poincare Conjecture}
\em The statement that $S^n$ is a strong Borel manifold is
equivalent to the \emph{Poincar\'e Conjecture} that every manifold
which is homotopy equivalent to a sphere $S^n$ is homeomorphic to
$S^n$. This follows from the fact that there are exactly two
homotopy classes of self-homotopy equivalences $S^n \to S^n$ which
both have homeomorphisms as representatives.  In particular $S^n$
is a Borel manifold if and only if it is a strong Borel manifold.
\em
\end{remark}

\begin{problem}[Classification of Borel manifolds] \label{pro:Classification_of_Borel_manifolds}
Which manifolds are Borel manifolds?
\end{problem}

In the light of both the Borel Conjecture and the Poincar\'e Conjecture, it is natural to consider the class
of manifolds $M$, whose universal covering $\widetilde{M}$ is homotopy-equivalent to a wedge
of $k$-spheres $S^k$ for some $2 \le k \le \infty$. We call such a manifold a \emph{generalized topological space forms}.
If $k \not= \infty$, this condition is equivalent to saying that the reduced integral homology
vanishes except in dimension $k$ and it is a direct sum of copies of $\bbZ$ in dimension $k$.
If $k = \infty$, then this condition is equivalent to saying that $M$ is an aspherical manifold.
A simply-connected generalized topological space form is the same as a homotopy sphere.
More generally, a generalized topological space form with finite fundamental group,
is the same as a  spherical topological space form.
If $G$ acts freely and cocompactly and properly discontinuously on $S^k \times \mathbb{R}^{m-k}$,
then $M= S^k \times \mathbb{R}^{m-k}/G$ is a generalized topological space form.
If $M$ and $N$ are $m$-dimensional aspherical manifolds, then $M \# N$ is a generalized topological space form.
If $M$ is aspherical, then for each  $k$ the manifold $M \times S^k$ is a generalized space form,
or more generally, all $S^k$-bundles over $M$  with $k>1$ are generalized space forms.

Most results in this paper concern generalized  space forms $M$. One can try to attack the question
whether $M$ is Borel by computing its \emph{structure set} $S^{top} (M)$. It consists of equivalence classes
of orientation preserving homotopy equivalences $N \to M$ with a manifold $N$ as source,
where two such homotopy equivalences $f_0 \colon N_0 \to M$ and $f_1 \colon N_1 \to M$ are equivalent
if there exists a homeomorphism $g \colon N_0 \to N_1$ with $f_1 \circ g\simeq f_0$.
The group $\hoaut_{\pi}(M)$ of homotopy
classes of self equivalences inducing the identity on $\pi_1$ up to conjugation acts on this
set by composition. A manifold $M$ is strongly Borel if and only if $S^{top} (M)$ consists of one element.
 A manifold $M$ is Borel if and only if  $\hoaut_{\pi}(M)$ acts transitively on $S^{top} (M)$ .
In general it is very hard to compute the structure set. But if the Farrell-Jones Conjecture
for $\pi_1 (M)$ holds, then one often can do this. More precisely we mean  the \emph{Farrell-Jones Conjecture}
for $K$- and $L$-theory
for the group $G$. In all relevant cases $G$ will be torsionfree. Hence this phrase will mean
that $\Wh(G)$ and $\widetilde{K}_n(\bbZ G)$ vanish for $n \le 0$ and that the assembly map
$H_n(BG;\bfL) \to L_n(\bbZ G)$ is bijective for all $n \in \bbZ$, where
$\bfL$ is the (non-connective periodic) $L$-theory spectrum and
$L_n(\bbZ G)$ is the $n$-th quadratic  $L$-group of $\bbZ G$. (We can ignore the
decoration by the Rothenberg sequences and the assumption that
 $\Wh(G)$ and $\widetilde{K}_n(\bbZ G)$ vanish for $n \le 0$.) More information about the
Farrell-Jones Conjecture can be found for instance in \cite{Farrell-Jones(1993a)},
\cite{Kreck-Lueck(2005)} and \cite{Lueck-Reich(2005)}.

For example the Farrell-Jones Conjecture holds for $\mathbb Z$ and the fundamental group of
surfaces of genus $\ge 1$. Combining this with the construction of certain self-equivalences,
we obtain the following result concerning generalized topological space forms.

\begin{theorem}[Sphere bundles over surfaces]
\label{the:Sphere_bundles_over_surfaces}
Let $K$ be $S^1$ or a $2$-dimen\-sio\-nal manifold different from $S^2$.
Let $S^d \to E \to K$ be a fiber bundle over $K$
such that $E$ is orientable and $d \ge 3 $.

Then $E$ is a Borel manifold. It is a  strong Borel manifold if and only if $K = S^1$.
\end{theorem}

This gives examples of Borel manifolds in all dimensions $> 3$, which are neither homotopy spheres nor aspherical.

In dimension $3$ the existence
of Borel manifolds is related to the Poincar\'e Conjecture and 
to Thurston's Geometrization Conjecture. Results of Waldhausen and Turaev imply:

\begin{theorem}[Dimension  $3$]\label{the:Thurston_and_property_Borel}
Suppose that Thurston's Geometrization Conjecture for irreducible $3$-manifolds
with infinite fundamental group
and the $3$-dimensional Poin\-ca\-r\'e Conjecture are true. Then every $3$-manifold with torsionfree
fundamental group is a strong Borel manifold.
\end{theorem}

Using the Kurosh theorem and the prime decomposition of $3$-manifolds one can even show that
if the assumptions of this theorem are fulfilled then the
fundamental group determines the homeomorphism type, a close analogy between
surfaces and $3$-manifolds (although the latter case is of course much more complicated).

Recently Perelman has announced a proof of Thurston's Geometrization Conjecture but details are still checked by the experts.

Given the analogy between the classification of surfaces and $3$-manifolds with torsionfree fundamental group,
it is natural to study in analogy to sphere bundles over surfaces sphere bundles over $3$-manifolds. Our result here is:

\begin{theorem}[Sphere bundles over $3$-manifolds]
\label{the:Sphere_bundles_over_3_manifolds}
Let $K$ be an aspherical $3$-dimensional manifold.
Suppose that the Farrell-Jones Conjecture for $K$- and $L$-theory holds for
$\pi_1(K)$. Let $S^d \to E \xrightarrow{\cong} K$ be a fiber bundle
over $K$ with orientable $E$
such that $d \ge 4$ or such that $d = 2,3$ and there is a map $i \colon
K \to E$ with $p \circ i \simeq \id_K$. Then

\begin{enumerate}

\item \label{the:Sphere_bundles_over_3_manifolds:i)}
$E$ is strongly Borel if and only if $H_1(K;\bbZ/2) = 0$;
\item \label{the:Sphere_bundles_over_3_manifolds:ii)}
If $d = 3\mod 4$ and $d \ge 7$, then $K \times S^d$ is Borel;
\item \label{the:Sphere_bundles_over_3_manifolds:iii)}
If $d = 0\mod 4$ and $d \ge 8$ and $H_1(K;\bbZ/2)\ne 0$, then $K \times
S^d$ is not Borel.
\end{enumerate}

\end{theorem}



The following result shows that if the fundamental groups of two $d$-dimen\-sional Borel manifolds $M$
and $N$ contain no $2-$torsion, then the connected sum $M \# N$ is a Borel manifold. Here we assume that $d> 4$.

\begin{theorem}[Connected sums]\label{the:connected_sums}
Let $M$ and $N$ be manifolds of the same dimension $n \ge 5$
such that neither $\pi_1(M)$ nor $\pi_1(N)$ contains elements of order $2$.
If both $M$ and $N$ are (strongly) Borel, then the same is true for their connected sum $M \# N$.

\end{theorem}

\begin{remark}
If $M$ and $N$ are aspherical Borel manifolds of dimension $\not= 4$ then $M \# N$ is a generalized topological space form, which is Borel.
\end{remark}

Combining the previous results, we have found infinitely many non-aspherical
and non-simply connected Borel manifolds in each dimension $\not= 4$.
The proof is in all cases based
on a determination of the structure set and by providing enough self equivalences following the scheme described above.

The main reason why these proofs do not work at present in dimension $4$ is that for the fundamental
groups under consideration it is not known whether they are good in the sense of Freedman. For this
reason one has to look at other classes of $4-$manifolds where also the determination of the structure
set is known but it is not clear how to construct enough self equivalences to guarantee a transitive action.
However,  if $\pi_1(M)$ is cyclic and $M$ is a spin manifold, one can use other methods to show that the homotopy
type determines the homeomorphism type (respecting the identification of fundamental groups). Previously
known Borel $4$ manifolds are flat $4$-manifolds, where the Borel Conjecture was proven.

\begin{theorem}[Dimension $4$] \label{the:dimension_4}
\begin{enumerate}

\item \label{the:dimension_4:finite_cyclic}
Let $M$ be a $4$-manifold with Spin structure such that its fundamental group is finite cyclic.
Then $M$ is Borel.

If $M$ is simply connected and Borel, then it has a Spin structure.

\item \label{the:dimension_4:flat} Let $N$ be a flat smooth Riemannian $4$-manifold or be $S^1 \times S^3$. Then
$N$ is strongly Borel. If $M$ is a simply connected $4$-manifold
with Spin-structure, then $M \# N$ is Borel.

\end{enumerate}
\end{theorem}

By Theorem~\ref{the:dimension_4}\eqref{the:dimension_4:flat}
we have provided infinitely
many non-aspherical and non-simply connected Borel manifolds in dimension $4$.
Except for $S^1 \times S^3$ and the flat manifolds, these manifolds are not generalized topological space forms.

 We have seen that under some mild restrictions the connected sum of two Borel
manifolds is a Borel manifold. It is natural to ask the corresponding question
for the cartesian product of Borel manifolds. If $M$ and $N$ are aspherical, then
$M \times N$ is aspherical and so Borel, if the Borel Conjecture holds. But if
the manifolds are not aspherical Borel manifolds, the picture becomes complicated.
An interesting test case is provided by the product of two spheres, where we give
a complete answer in terms of the unstable Arf invariant.

Let $\Omega_{k,k+d}^{\fr}$ be the bordism set of smooth $k$-dimensional manifolds
$M$ with an embedding $i \colon M \to \bbR^{k+d}$ together with an  (unstable) framing
of the normal bundle $\nu(i)$. If $d > k$, this is the same as the
bordism group $\Omega_k^{\fr}$ of stably framed smooth $k$-dimensional manifolds
since any $k$-dimensional smooth manifold $M$ admits an embedding into $\bbR^{k+d}$
and a stable framing on $\nu(M,\bbR^{k+d})$ is the same as an unstable framing for $d > k$.
The \emph{Arf invariant} $\Arf(M) \in \bbZ/2$ of a stably framed manifold $M$ is the Arf invariant of the surgery problem associated
to any degree one map $M \to S^{\dim(M)}$ with the obvious bundle data coming from the stable framing. It induces a homomorphism
of abelian groups
\begin{eqnarray}
\Arf_k \colon \Omega_k^{\fr} & \to & \bbZ/2.
\label{stable_Arf_homomorphism}
\end{eqnarray}
If $g \colon S^k \times S^d \to S^k$ is a map, define its
\emph{Arf invariant} $\Arf(g) \in \bbZ/2$ to be the Arf invariant
of the stably framed manifold $\overline{g}^{-1}(\pt)$ for any map
$\overline{g} \colon S^k \times S^d \to S^k$ which is homotopic to
$g$ and transverse to $\pt \subseteq S^d$. Here the stable
framing of $\overline{g}^{-1}(\pt)$ is given by  the standard
stable framing of the normal bundle of $S^k \times S^d$ and the
trivialization of the normal bundle of $\overline{g}^{-1}(\pt)
\subseteq S^k \times S^d$ coming from transversality.

\begin{theorem}[Products of two spheres]

\label{the:S^k_times_S^d}
Consider $k,d \in \bbZ$ with $k,d \ge 1$. Then:
\begin{enumerate}

\item \label{the:S^k_times_S^d:strong_Borel}
Suppose that $k +d \not= 3$. Then $S^k \times S^d$ is a strong Borel manifold if and only if both
$k$ and $d$ are odd;

\item \label{the:S^k_times_S^d:S^1timesS^2} For $d \ne 2$ the manifolds $S^1\times S^d$ is Borel, and
$S^1 \times S^2$ is strongly Borel if and only if the $3$-dimensional Poincar\'e Conjecture is true;

\item \label{the:S^k_times_S^d:S^2timesS^2}
The manifold $S^2 \times S^2$ is Borel but not strongly Borel;

\item  \label{the:S^k_times_S^d:necessary_and_sufficient_for_Borel}
Suppose $k,d >1$ and $k +d \ge 4$. Then the manifold $S^k \times S^d$ is Borel if and only if the following conditions are satisfied:

\begin{enumerate}

\item \label{the:S^k_times_S^d:necessary_and_sufficient_for_Borel:values_for_k_and_d}
Neither $k$ nor $d$ is divisible by $4$;

\item \label{the:S^k_times_S^d:necessary_for_Borel:Arf_invariant_for_k}
If $k = 2 \mod 4$, then there is a map $g_k \colon S^k \times S^d \to S^k$
such that its Arf invariant $\Arf(g_k)$ is non-trivial and
its restriction to $S^k \times \pt$ is an orientation preserving
homotopy equivalence $S^k \times \pt \to S^k$;

\item \label{the:S^k_times_S^d:necessary_for_Borel:Arf_invariant_for_d}
If  $d = 2 \mod 4$, then there is a map $g_d \colon S^k \times S^d \to S^d$
such that its Arf invariant $\Arf(g_d)$ is non-trivial and
its restriction to $ \pt \times S^d$ is an orientation preserving
homotopy equivalence $\pt \times S^d\to S^d$.

\end{enumerate}

\end{enumerate}

\end{theorem}

\begin{remark}[Relation to the Arf-invariant-one-problem]
\label{rem:Relation_to_the_Arf-invariant-one-problem} \em
The con\-di\-tion~\eqref{the:S^k_times_S^d:necessary_for_Borel:Arf_invariant_for_k}
appearing in
Theorem~\ref{the:S^k_times_S^d}~\eqref{the:S^k_times_S^d:necessary_and_sufficient_for_Borel}
implies that the (stable) Arf invariant homomorphism $\Arf_k$ of
\eqref{stable_Arf_homomorphism} is surjective. The problem whether
$\Arf_k$ is surjective is the famous \emph{Arf-invariant-one-problem}
(see \cite{Browder(1969)}). The map $\Omega_k^{\fr} \to \bbZ/2$ is
known to be trivial unless $2k+2$ is of the shape $2^l$ for some
$l \in \bbZ$ (see \cite {Browder(1969)}). Hence a necessary
condition for $S^k \times S^d$ to be Borel is that $k$ is odd or
that $2k+2$ is of the shape $2^l$ for some $l \in \bbZ$ and
analogously for $d$.

Suppose that the unstable Arf-invariant-map
\begin{eqnarray}
\Arf_{k,k+d} \colon \Omega_{k,k+d}^{\fr} & \to & \bbZ/2.
\label{unstable_Arf_homomorphism}
\end{eqnarray}
is surjective. Then
condition~\eqref{the:S^k_times_S^d:necessary_for_Borel:Arf_invariant_for_k}
is automatically satisfied by the following argument. Choose a
framed manifold $M \subseteq S^{k+d}$ with $\Arf_{k,k+d}([M]) =
1$. By the Pontrjagin-Thom construction we obtain a map $g_k'
\colon S^{k+d} \to S^k$ which is transversal to $\pt \subseteq
S^k$ and satisfies $g^{-1}(\pt) = M$. Now define the desired map
$g_k$ to be the composition of $g_k'$ with an appropriate collapse
map $S^k \times S^d \to S^{k+d}$. So the surjectivity of the
unstable Arf invariant map~\eqref{unstable_Arf_homomorphism}
implies
condition~\eqref{the:S^k_times_S^d:necessary_for_Borel:Arf_invariant_for_k}.
Of course the surjectivity of~\eqref{unstable_Arf_homomorphism} is
in general a stronger condition than the surjectivity
of~\eqref{stable_Arf_homomorphism}. The unstable Arf invariant
map~\eqref{unstable_Arf_homomorphism} is known to be surjective if
$k=2$ and $d \ge 1$. Hence $S^k \times S^d$ is Borel if $k = 2$
and $d \ge 2$ or if $k \ge 2$ and $d = 2$.
\end{remark}

Now we discuss the following question. How complicated can the homotopy type of Borel manifolds be?
In the situation of the Borel and Poincar\'e Conjectures the homotopy type is determined by the
fundamental group and - in the case of homotopy spheres - by the homology groups. Most of our
results concerned generalized topological space forms (or connected sums of these), spaces
 which are ''close neighbours'' of aspherical manifolds resp. homotopy spheres. Besides the
products of spheres we have given results concerning other classes of manifolds only in
dimension $4$. If we concentrate on manifolds with torsionfree fundamental groups
(the lens spaces show that even very simple manifolds are in general not Borel for
cyclic fundamental groups (see for instance \cite[\S~29 and \S~31]{Cohen(1973)},  \cite[Section 2.4]{Lueck(2002c)})),
these results in dimension $4$ for non-aspherical manifolds concern manifolds with fundamental group $\mathbb {Z}$.
Here the fundamental group and the intersection form on $\pi_2$, which is a homotopy invariant, determines
the homeomorphism type. The following classes of manifolds are natural generalizations of these manifolds.

\begin{problem}[Classification of certain low-dimensional manifolds]
\label{pro:Classification_of_certain_low-dimensional_manifolds}
Classify up to orientation preserving homotopy equivalence,
homeomorphism (or diffeomorphism in the smooth case) all
manifolds in dimension $1 \le k <  n \le 6$ for which $\pi = \pi_1(M)$
is non-trivial and is isomorphic to $\pi_1(K)$ for a manifold $K$ of dimension $k \le 2$ with $\pi_1(K) \not= \{1\}$
and whose second homotopy group $\pi_2(M)$ vanishes.
\end{problem}

\begin{remark}[Simply-connected case] \label{rem:simply-connected_case} \em
We have excluded in Problem~\ref{pro:Classification_of_certain_low-dimensional_manifolds} the case
$\pi_1(K) = \{1\}$ since then a complete answer to this problem is already known.
Namely, if $M$ is a $2$-connected $n$-dimensional manifold, then $n \ge 3$,
it is homotopy equivalent to $S^3$ for $n = 3$ and it is homeomorphic to $S^n$ for $n = 4,5$.
If $n = 6$ and $M$ is a $2$-connected smooth manifold, then its oriented homotopy type and its oriented
diffeomorphism  are determined by the intersection from on $H_3(M)$
(see Wall~\cite{Wall(1962)}). This also applies to the topological category
by the work of Kirby-Siebenmann~\cite{Kirby-Siebenmann(1977)}.
\em
\end{remark}

The following results give an almost complete answer to this problem.

\begin{theorem}{\bf (Manifolds appearing in Problem~\ref{pro:Classification_of_certain_low-dimensional_manifolds} of dimension $\le 5$).}
\label{the:motivating_problem_n_le_5}
Let $M$ and $K$ be as in
Problem~\ref{pro:Classification_of_certain_low-dimensional_manifolds}.
Let $f \colon M \to K = B\pi$ be the classifying map
for the $\pi$-covering $\widetilde{M} \to M$.
Suppose $n \le 5$. Then

\begin{enumerate}

\item \label{the:motivating_problem_n_le_5:homotopy_classification}
Both the oriented homotopy type and the oriented homeomorphism type of $M$ are determined by its second Stiefel-Whitney class.
Namely, there is precisely one fiber bundle $S^3 \to E \xrightarrow{p} K$ with structure group $SO(4)$ whose
second Stiefel-Whitney class agrees with the second Stiefel-Whitney class of $M$
under the isomorphism $H^2(f;\bbZ/2) \colon H^2(K;\bbZ/2) \xrightarrow{\cong}
H^2(M;\bbZ/2)$. There exists an orientation preserving homeomorphism $g \colon M \to E$
such that $p \circ g$ and $f$ are homotopic.

In particular the vanishing of the second Stiefel-Whitney class of
$M$ implies that there is an orientation preserving homomorphism
$g \colon M \to S^3 \times B\pi$ such that $\pr_{B\pi} \circ g$
and $f$ are homotopic;

\item \label{the:motivating_problem_n_le_5:Borel}
The manifold $M$ is never a strong Borel manifold but is always a Borel manifold.

\end{enumerate}
\end{theorem}

\begin{theorem}{\bf (Manifolds appearing in Problem~\ref{pro:Classification_of_certain_low-dimensional_manifolds} of dimension $6$).}
\label{the:motivating_problem_n=6}
Let $M$ and $K$ be as in
Problem~\ref{pro:Classification_of_certain_low-dimensional_manifolds}.
Suppose that $n = 6$.

\begin{enumerate}
\item \label{the:motivating_problem_n=6:Borel}
The manifold $M$ is never strongly Borel but always Borel;

\item \label{the:motivating_problem_n=6:Spin_case}
Suppose that $w_2(M) = 0$, or, equivalently, that $M$ admits a Spin-structure.
Then both the oriented homotopy type and oriented homeomorphism type of $M$ is determined
by the $\bbZ\pi$-isomorphism class of the intersection form
on $H_3(\widetilde{M})$.

Given a finitely generated stably free $\bbZ\pi$-module
together with a non-degenerate skew-symmetric $\bbZ\pi$-form on it,
it can be realized as the intersection form of a $6$-dimensional manifold having the properties described in
Problem~\ref{pro:Classification_of_certain_low-dimensional_manifolds}.

\end{enumerate}
\end{theorem}

One expects that Borel manifolds are an exception. The following results which give necessary conditions for $M$ to be Borel
support this intuition.

\begin{theorem}[A necessary condition for sphere bundles over aspherical manifolds] \label{the:necessary_condition_for_sphere bundles}
Let $S^d \to E \xrightarrow{p} K$ be a fiber bundle such that $E$ and $K$ are manifolds and
$K$ is aspherical. Assume that there is a map $i \colon K \to E$ with $p \circ i \simeq \id_K$.
Suppose that $d$ is odd and $d \ge 3$ or suppose that $d$ is even and $k \le d-2$ for $k = \dim(K)$.

Then a necessary condition for $E$ to be a Borel manifold is that
$H_{k-4i}(K;\bbQ)$ vanishes for all $i \in \bbZ, i \ge 1$.
\end{theorem}

\begin{theorem}[A necessary condition for being a Borel manifold]
\label{the:A_necessary_condition_for_being_a_Borel_manifold}
Let $M$ be a Borel manifold of dimension $n$ with fundamental group $\pi$. Let $\hoaut_{\pi}(M)$ be the set
of homotopy classes of orientation preserving self-ho\-mo\-to\-py equivalences
$f \colon M \to M$ which induce up to conjugation the identity on the fundamental group.
Let $\call(M)_{i} \in H_{4i}(M;\bbQ)$ be the $i$-th component
of the $L$-class $\call(M)$ of $M$.

Then  the subset of $\bigoplus_{i \in \bbZ, i \ge 1} H_{n-4i}(M;\bbQ)$
$$S: = \{f_*(\call(M) \cap [M]) - \call(M) \cap [M] \mid [f] \in \hoaut_{\pi}(M)\}$$
is an abelian subgroup and the $\bbQ$-submodule generated by $S$
must contain the kernel of the map induced by the classifying map $c \colon M \to B\pi$
$$c_* \colon \bigoplus_{i \in \bbZ, i \ge 1} H_{n-4i}(M;\bbQ) \to
\bigoplus_{i \in \bbZ, i \ge 1} H_{n-4i}(B\pi;\bbQ) .$$

In particular for every $i \ge 1$ with $\call(M)_{i} = 0$ the map
$c_* \colon H_{n-4i}(M;\bbQ) \to H_{n-4i}(B\pi;\bbQ)$ must be injective.

\end{theorem}

We have mentioned that lens spaces are in general not Borel, an  indication that torsion
in the fundamental group makes Borel less likely. The following result shows that in dimension $4k + 3$ torsion excludes Borel.

\begin{theorem}[Chang-Weinberger\cite{Chang-Weinberger(2003b)}]
\label{intr: Chang-Weinberger}
Let $M^{4k+3}$ be a closed oriented manifold for $k \ge 1$ whose fundamental
group has torsion. Then there are  infinitely many pairwise not homeomorphic smooth manifolds
which are homotopy equivalent to M (and even simply and tangentially homotopy
equivalent to $M$) but not homeomorphic to M.
\end{theorem}

Another natural class of manifolds are the homology spheres where surgery gives a necessary and sufficient condition for being strongly Borel:

\begin{theorem}[Homology spheres]
\label{the:Integral_homology_spheres}
Let $M$ be a $n$-dimensional manifold of dimension $n \ge 5$ with
fundamental group $\pi = \pi_1(M)$.

\begin{enumerate}

\item \label{the:Integral_homology_spheres:L-criterion}

Let $M$ be an integral homology sphere.
Then $M$ is strongly Borel if and only if
the inclusion $j \colon \bbZ \to \bbZ \pi$
induces an isomorphism
$$L_{n+1}^s(j) \colon L_{n+1}^s(\bbZ) \xrightarrow{\cong} L_{n+1}^s(\bbZ \pi)$$
on the simple $L$-groups $L^s_{n+1}$;

\item \label{the:Integral_homology_spheres:nec_n+1}
Suppose that $M$ is a rational homology sphere and Borel.
Suppose that $\pi$ satisfies the Novikov Conjecture in dimension $(n+1)$, i.e. the assembly map
$H_{n+1}(B\pi;\bfL) \to L_{n+1}(\bbZ \pi)$ is rationally injective. Then
$$H_{n+1-4i}(B\pi;\bbQ) = 0$$
for $i \ge 1$ and $n+1-4i \not= 0$.
\end{enumerate}
\end{theorem}

The paper is organized as follows:\\[1mm]
\begin{tabular}{ll}
\ref{sec:On_the_Structure_Set_of_Certain_Topological_Manifolds}.
&
On the Structure Set of Certain Topological Manifolds
\\
\ref{sec:Constructing_Self-homotopy_Equivalences}.
&
Constructing Self-homotopy Equivalences
\\
\ref{sec:Sphere_Bundles}
&
Sphere Bundles
\\
\ref{sec:connected_sums}
&
Connected Sums
\\
\ref{sec:Dimension_3}
&
Dimension $3$
\\
\ref{sec:Dimension_4}
&
Dimension $4$
\\
\ref{sec:Products_of_Two_Spheres}.
&
Products of Two Spheres
\\
\ref{sec:On_the_Homotopy_Type_of_Certain_Low-Dimensional_Manifolds}.
&
On the Homotopy Type of Certain Low-Di\-men\-sio\-nal Manifolds
\\
\ref{sec:On_the_Classification_of_Certain_Low-D-mensional_Manifolds}.
&
On the Classification of Certain Low-Di\-men\-sio\-nal Manifolds
\\
\ref{sec:A_Necessary_Condition_for_Being_a_Borel_Manifold}.
&
A Necessary Condition for Being a Borel Manifold
\\
\ref{sec:Integral_Homology_Spheres} &
Integral Homology Spheres
\\
&
References
\end{tabular}
\\[2mm]
We thank Andrew Ranicki for fruitful discussions about this paper.


\typeout{--------------------------------   Section 1 ------------------------------------}

\section{On the Structure Set of Certain Topological Manifolds}
\label{sec:On_the_Structure_Set_of_Certain_Topological_Manifolds}

We begin with a fundamental criterion for Borel manifolds
which follows directly from the definitions.

\begin{theorem}[Surgery criterion for Borel manifolds]\label{the:surgery_criterion_for_Borel}
A manifold $M$ is a Borel manifold if and  only if the action of the group of homotopy classes
of self-homotopy equivalence $M \to M$ which induce the identity on the fundamental group up to conjugation
on the topological structure set $\cals^{\topo}(M)$ is transitive, and $M$ is a strong
Borel manifold if and only $\cals^{\topo}(M)$ consists of one element.
\end{theorem}

Now we determine the topological structure set for certain manifolds.
In the sequel we denote by $\bfL\langle 1 \rangle$ the $1$-connected cover of the
quadratic $L$-theory spectrum $\bfL$ and by $\bfu \colon \bfL\langle 1 \rangle \to \bfL$ the canonical map.
We get
\begin{multline*}
\pi_q(\bfL\langle 1 \rangle)
~ = ~
\left\{\begin{array}{ll}
L_q(\bbZ) &  \text{, if } q \ge 1\\
0         &  \text{, if } q \le 0
\end{array}\right\}
~ = ~
\left\{\begin{array}{ll}
\bbZ   &  \text{, if } q \ge 4, q = 0 \mod 4\\
\bbZ/2 &   \text{, if } q \ge 2, q = 2 \mod 4\\
0      &   \text{ otherwise}
\end{array} \right.
\end{multline*}

\begin{theorem} \label{the:computing_the_structure_set}
Let $M$ be an $n$-dimensional manifold for $n \ge 5$.
Let $K$ be an  aspherical $k$-dimensional manifold with fundamental group $\pi$.
Suppose that the Farrell-Jones Conjecture for algebraic $K$- and $L$-theory
holds for $\pi$. Let $f \colon M \to K$ be a $2$-connected map.
Suppose that we can choose a  map $i \colon K \to M$
such that $f \circ i$ is homotopic to the identity.

\begin{enumerate}

\item \label{the:computing_the_structure_set:H_m(K;bfu)}
The homomorphism
$$H_m(\id_K;\bfu) \colon H_m(K;\bfL\langle 1 \rangle) \to H_m(K;\bfL)$$
is bijective for $m \ge k+1$ and injective for $k = m$;

\item \label{the:computing_the_structure_set:structure_set}
The exact topological surgery sequence
for $M$ yields the short split-exact sequence
$$0 \to  \cals^{\topo}(M) \xrightarrow{\sigma_n}
H_n(M;\bfL\langle 1 \rangle) \xrightarrow{H_n(f;\bfL\langle 1 \rangle)}
H_n(K;\bfL\langle 1 \rangle) \to 0.$$
In particular we get an isomorphism
$$\overline{\sigma}_n \colon \cals^{\topo}(M) \xrightarrow{\cong}
H_n(i \colon K \to M;\bfL\langle 1 \rangle).$$
\end{enumerate}
\end{theorem}
\begin{proof}
\eqref{the:computing_the_structure_set:H_m(K;bfu)}
Let $\bfE$ be the homotopy fiber of $\bfu$. Hence we have a fibration of spectra
$$\bfE \to \bfL\langle 1 \rangle \xrightarrow{\bfu} \bfL$$
which induces a long exact sequence
\begin{multline*}
\ldots \to H_{m+1}(K;\bfL\langle 1 \rangle) \to H_{m+1}(K;\bfL) \to H_m(K;\bfE)
\\
\to H_m(K;\bfL\langle 1 \rangle) \to H_m(K;\bfL) \to \ldots.
\end{multline*}
Since $\pi_q(\bfE) = 0$ for $q \ge 0$, an easy spectral sequence argument shows that
$H_m(K;\bfE) = 0$ for $m \ge k$. Hence the map
$$H_m(\id_K;\bfu) \colon H_m(K;\bfL\langle 1 \rangle) \to H_m(K;\bfL)$$
is bijective for $m \ge k+1$ and injective for $k = m$.
\\[1mm]
\eqref{the:computing_the_structure_set:structure_set}
There is an exact sequence of abelian groups called algebraic surgery exact sequence
\cite[Definition~15.19 on page~169]{Ranicki(1992)}.
\begin{multline}
\ldots  \xrightarrow{\sigma_{n+1}}  H_{n+1}(X;\bfL\langle 1 \rangle)
\xrightarrow{A_{n+1}} L_{n+1}(\bbZ\pi_1(X)) \xrightarrow{\partial_{n+1}}
\\
 \cals^{\topo}(X) \xrightarrow{\sigma_n}
H_{n}(X;\bfL\langle 1 \rangle) \xrightarrow{A_{n}} L_{n}(\bbZ\pi_1(X)) \xrightarrow{\partial_{n}}  \ldots
\label{Ranickis algebraic surgery sequence}
\end{multline}
which is  defined for every simplicial connected complex $X$ and natural in $X$.
It  agrees for $X$ a $n$-dimensional manifold for $n \ge 5$ with the Sullivan-Wall geometric
exact surgery sequence \cite[Theorem~18.5 on page~198]{Ranicki(1992)}.
Notice that by assumption
$\Wh_i(\pi_1(K))$ vanishes for $i \le 1$ so that we can ignore any decorations
in the sequel.
The following diagram commutes for all $m$
$$
\begin{CD}
H_m(M;\bfL\langle 1 \rangle) @> A_m >>  L_m(\bbZ\pi)
\\
@V H_m(f;\bfL\langle 1 \rangle)  VV @V \id V \cong V
\\
H_m(K;\bfL\langle 1 \rangle) @> A_m >>  L_m(\bbZ\pi)
\\
@V H_m(\id_K; \bfu)  VV @V \id V \cong V
\\
H_m(K;\bfL) @> A_m >\cong >  L_m(\bbZ\pi)
\end{CD}
$$
The existence of the map $i$ with $f \circ i \simeq \id_K$ ensures
that $H_m(f;\bfL\langle 1 \rangle)$ is surjective.
Now the claim follows from
assertion~\eqref{the:computing_the_structure_set:H_m(K;bfu)} and the exact sequence
\eqref{Ranickis algebraic surgery sequence} for $X = M$.
\end{proof}

\begin{theorem} \label{the:structure_set_of_S^d-bundles_over_K}
Let $K$ be an  aspherical $k$-dimensional manifold with fundamental group $\pi$.
Suppose that the Farrell-Jones Conjecture for algebraic $K$- and $L$-theory
holds for $\pi$. Let $d$ be an integer with $d \ge 2$ and $d + k \ge 5$.
Consider a fiber bundle $S^d \to E \xrightarrow{p} K$ such that $E$ is oriented.
Suppose that there exists a map $i \colon K \to E$ with $p \circ i = \id$.

Then $i \colon K \to E$ is an embedding of topological manifolds and
we obtain an isomorphism
of abelian groups
$$\cals^{\topo}(E) \xrightarrow{\cong}
H_k(K;\bfL\langle 1 \rangle).$$ It sends under the identification
of  $H_k(K;\bfL\langle 1 \rangle)$ with the set $\caln(K) \cong
[K,G/TOP]$ of normal surgery problems with $K$ as reference space
an element $f \colon M \to E$ to the following surgery problem: By
changing $f$ up to homotopy we can arrange that $f$ is transverse
to $i \colon K \to E$. Let $g \colon N = f^{-1}(i(K)) \to K$ be
the map of manifolds of degree $1$
induced by $f$ and $i^{-1} \colon i(K) \to K$. By transversality
we obtain a bundle map $\overline{g} \colon \nu(N,M) \to \nu(i)$
covering $g$. Choose a vector bundle $\xi \to M$ and a bundle map
$\overline{f} \colon \nu(M) \to \xi$ covering $f \colon M \to E$.
Then $g$ is covered by the bundle map
$$\overline{g} \oplus \overline{f}|_N \colon \nu(N) = \nu(N,M) \oplus \nu(M)|_N \to \nu(i) \oplus i^*\xi,$$
and these data give the desired surgery problem with target $K$.
\end{theorem}

\proof See \cite[Chapter~11]{Wall(1999)}, \cite[pages~257--260]{Ranicki(1992)}. \qed

\begin{theorem} \label{the:structure_set_for_S^d_bundle_over_K_with_dim(K)le3}
Let $K$ be an  aspherical $k$-dimensional manifold with fundamental group $\pi$.
Consider a fiber bundle $S^d \to E \xrightarrow{p} K$ with $d \ge 1$.
Suppose that the Farrell-Jones Conjecture for algebraic $K$- and $L$-theory
holds for $\pi$.
Then:

\begin{enumerate}

\item \label{the:structure_set_for_S^d_bundle_over_K_with_dim(K)le3:k=2}

If $k = 2$ and  $d \ge 3$, then
$$\cals^{\topo}(E) \cong L_2(\bbZ) \cong \bbZ/2;$$

\item \label{the:structure_set_for_S^d_bundle_over_K_with_dim(K)le3:k=3}
If $k = 3$ and $d \ge 4$ or if $k = 3$, $d = 2,3$ and there is a map $i \colon K \to E$ with $p \circ i = \id_K$,
then $$\cals^{\topo}(E) \cong H_1(K;L_2(\bbZ)) \cong H_1(K;\bbZ/2);$$

\item \label{the:structure_set_for_S^d_bundle_over_K_with_dim(K)le3:k=1}
If $k=1$ and $d \ge 3$, then
$$\cals^{\topo}(E) = 0;$$

\end{enumerate}
\end{theorem}
\begin{proof}
We first prove the claim in the case $k+d \ge 5$.
Since $K$ is $k$-dimensional and $S^d$ is $(d-1)$-connected,
we can find a map $i \colon K \to E$ such that $p \circ i$ is homotopic to $\id_K$ provided
$d > k$. By the homotopy lifting property we can arrange that $p \circ i$ is $\id_K$.
By assumption such a map $i$ exists also in the remaining case $k = 3$ and $d = 2,3$.
Now the claim follows from Theorem~\ref{the:structure_set_of_S^d-bundles_over_K} and
an easy computation with the Atiyah-Hirzebruch spectral sequence.
Thus we have proven assertions
\eqref{the:structure_set_for_S^d_bundle_over_K_with_dim(K)le3:k=2} and
\eqref{the:structure_set_for_S^d_bundle_over_K_with_dim(K)le3:k=3} and
for assertion \eqref{the:structure_set_for_S^d_bundle_over_K_with_dim(K)le3:k=1}
only the case $k = 1$ and $d = 3$ remains open.
Then $\pi_1(E)$ is $\bbZ$ which is a good fundamental group in the sense
of Freedman~\cite{Freedman(1983)}.
Hence topological surgery works also in this dimension $4$ and the same argument
which gives the claim for $k = 1$ and $k + d \ge 5$ works also for
$k = 1$ and $k +d = 4$.
\end{proof}


\typeout{--------------------------------   Section 2 ------------------------------------}

\section{Constructing Self-homotopy Equivalences}
\label{sec:Constructing_Self-homotopy_Equivalences}

In this section we describe a certain construction of self-homotopy equivalences.
It will be used to show that the action of the group of orientation preserving
homotopy equivalences $E \to E$ which induce the identity on the fundamental group up to conjugation
on $\cals^{\topo}(E)$ is transitive for certain manifolds $E$.

Suppose we are given the following data:

\begin{itemize}

\item Let $K$ and $E$ be manifolds with $\dim(K) = k$ and $\dim(E) = k+d$ for $k,d \ge 2$;

\item An embedding $i_K \colon K \to E$;

\item A map $\phi \colon S^d \to E$ which is transversal to $i_K \colon K \to E$ and
the intersection of the images $\im(\phi)$ and $\im(i_K)$ consists of precisely one point $e_0 \in E$;

\item Let $M$ be a manifold together with an embedding
$i_M \colon M \to S^{k+d}$ and a framing of the normal bundle $\mu(i_M)$;

\item $k +d \ge 5$ or $E$ is simply connected.

\end{itemize}

Fix an embedding $j_K \colon D^k \to K$ with $i_K \circ j_K(0) = e_0$.
Since $D^k$ is contractible, we can choose a disk bundle map
\comsquare{D^k \times D^d}{\overline{j_K}}{D\nu(i_K)}{\pr}{}
{D^k}{j_K}{K}
which is fiberwise a homeomorphism. Choose an embedding
$j_M \colon D^k \to M$. Since $D^k$ is contractible, we can choose
a disk bundle map
\comsquare{D^k \times D^d}{\overline{j_M}}{D\nu(M,S^{k+d})}{\pr}{}
{D^k}{j_M}{M}
which is fiberwise a homeomorphism. Using a tubular neighborhood, we will also regard
$\overline{j_K}$ as an embedding $D^k \times D^d \to E$ extending
$j_K \colon D^k \to E$ and
$\overline{j_M}$ as an embedding $D^k \times D^d \to S^{k+d}$ extending
$j_M \colon D^k \to S^{k+d}$.

In the sequel let $j_{S^d} \colon  D^d \to S^d$ be the obvious embedding given by the lower hemisphere.
Recall that $\phi$ and $i_K$ are transversal to one another and the intersection of their
images consists of the point $e_0$.
We can assume without loss of generality that $\phi \circ j_{S^d}(0) = e_0 = j_K(0)$ holds.
Now we can thicken $\phi$ to  a map $\overline{\phi} \colon D^k \times S^d \to E$
such that the composite
$$D^k \times D^d \xrightarrow{\id \times j_{S^d}} D^k \times S^d \xrightarrow{\overline{\phi}} E$$
agrees with the embedding $\overline{j_K} \colon D^k \times D^d \to E$ and the intersection of
the image of $\overline{\phi}$ and of $D\nu(i_K)$ considered as subset of $E$ is
the image of $\overline{j_K} \colon D^k \times D^d \to E$.

The Pontrjagin-Thom construction applied to $i_M \colon M \to
S^{d-k}$ together with the given framing on $M$ yields a map $\PT
\colon S^{k+d} \to S^d$ such that $\PT$ is transversal to
$i_{S^d}(0) \in S^d$, the preimage of $i_{S^d}(0) \in S^d$ is just
$M$ and the bundle map given by transversality $\nu(i_M) \to
\nu(\{i_{S^d}(0)\} \subseteq S^d)$ is just the given framing. We
can arrange by shrinking $D^d$ that the composition of $\PT$ with
$\overline{j_M} \colon D^k \times D^d \to S^{k+d}$ agrees with the
composite
$$D^k \times D^d \xrightarrow{\pr} D^d \xrightarrow{j_{S^d}} S^d.$$
Choose a map $c \colon S^{k+d} \to D^k$ such that its composite
with the embedding $\overline{j_M} \colon D^k \times D^d \to
S^{k+d}$ is the projection $\pr \colon D^k \times D^d \to D^k$, it
is transversal to $0 \in D^k$ and the preimage of $0$ is an
embedded $S^d$. Such a map can be constructed from the map
$S^{k+d} = S^d \ast S^{k-1} \to \pt \ast S^{k-1} = D^k$.

Now consider the composite
$$\alpha \colon S^{k+d} \xrightarrow{c \times \PT} D^k \times S^d \xrightarrow{\overline{\phi}} E.$$
Its composition with the embedding $\overline{j_M} \colon D^k
\times D^d \to S^{k+d}$ agrees with the embedding $\overline{j_K}
\colon D^k \times D^d \to E$. It is transverse to $i_K \colon K
\to E$ such that the preimage of both $i_K(K)$ and of $i_K \circ
j_K(D^k)$ is $M$ and the bundle map given by transversality from
$\nu(M,S^{k+d}) \to i_K^*\nu(i)$ is compatible up to isotopy with
the given framing of $\nu(M;S^{k+d})$ and some framing of the
bundle $i_K^*\nu(i)$ over $D^k$.

In the sequel we consider the connected sum $E \# S^{n+k}$ with respect to the two embeddings
$\overline{j_M} \colon D^k \times D^d \to S^{k+d}$ and
$\overline{j_K} \colon D^k \times D^d \to E$. By construction the identity
$\id \colon E \to E$ and the map $\alpha \colon S^{k+d} \to E$ fit together and yield
a map $\id \# \alpha \colon E \# S^{k+d} \to E$.

We claim that this map is a homotopy equivalence. Choose a point $x \in E$ which not contained in
the images of $\overline{j_K} \colon D^k \times D^d  \to E$ and of
$\overline{\phi} \colon D^k \times S^d \to E$. Then the preimage of $x$ under
$\id \# \alpha$ is $x$ and the map $\id \# \alpha$ induces the identity on a neighborhood
of $x$. This implies that $\id \# \alpha$  has degree one. The inclusions
of $E- \im(\overline{j_K})$ into both $E$ and $E \# S^{k+d}$ are $(k+d-1)$-connected.
Since $\id \# \alpha$ induces the identity on $E- \im(\overline{j_K})$,
the homomorphisms $\pi_j(\id \# \alpha)$ are bijective for $j \le k+d-2$. By assumption we have $k +d \ge 5$
or we have $k+d = 4$ and $\pi_1(E) = \{1\}$. Now we conclude from Poincar\'e duality that
$\id \# \alpha \colon E \# S^{k+d} \to E$ is a homotopy equivalence.

Obviously we can find a homeomorphism $\beta \colon E \to E \#
S^{k+d}$ such that the composite $\alpha \circ \beta$ is the identity
outside the image of $\overline{j_K} \colon D^k \times D^d \to E$.

The  map $\id \# \alpha \colon E \# S^{k+d} \to E$ is transversal
to $i_K \colon K \to E$. The preimage of $i_K$ is the connected
sum $K \# M$, which is taken  with respect to the embeddings $j_M
\colon D^k \to M$ and $j_K \colon D^k \to K$. This map has degree
one and is covered by bundle data due to transversality. The
resulting normal map with target $E$ agrees with the connected sum
of the normal map $\id \colon K \to K$ and the normal map $M \to
S^k$ given by the collapse map of degree one and the bundle data
coming from the given framing on $M$.

Now additionally suppose that the given map $\phi \colon S^d \to
E$ is an embedding. (It is automatically a local embedding near
the intersection point with $K$ by transversality but a priori not
a global embedding). Then  also the map $\overline{\phi} \colon
D^k \times S^d \to E$ can be chosen to be an embedding. It is not
hard to check that the map $\id \# \alpha \colon E \# S^{k+d} \to
E$ is transversal to $\phi$ and the corresponding surgery problem
is given by a homeomorphism $(\id \# \alpha)^{-1}(\alpha(S^d))
\xrightarrow{\cong} S^d$ covered by an isomorphism of the normal
bundles. In particular this surgery problem represents the trivial
element in $\caln(S^d)$.

Next we explain the maps in the following diagram
\begin{eqnarray}
& \comsquarewithout{\Omega^{\fr}_{k,k+d}}{\tau}{\cals^{\topo}(E)}
{a}{b}
{\caln(S^k)}{cs_K}{\caln(S^k)}
&
\label{diagram_for_aut_and_structure_set_for_K_subseteq_E}
\end{eqnarray}
Recall that we denote by $\Omega^{\fr}_{k,d}$ the set of bordism classes of
$k$-dimensional manifolds $M$ together with an embedding $M
\subseteq \bbR^{k+d}$ and an (unstable)  framing of its normal bundle
$\nu(M\subseteq \bbR^{k+d})$. The map $\tau$ is given by the
construction above which assigns to $[M] \in \Omega^{\fr}_{k,k+d}$
the element in the structure set given by $\alpha \colon E \#
S^{k+d} \to E$. The map $a$ sends a framed $k$-dimensional
manifold $M$ to the normal map given by the collapse map $c \colon
M \to S^k$ covered by stable bundle map from $\nu(M)$ to the
trivial bundle over $S^k$ given by the framing. The map $cs_K$ is
given by taking the connected sum of a surgery problem with target
$S^k$ with the one given by the identity $\id_K \colon K \to K$.
The map $b$ sends the class of a homotopy equivalence $f \colon M
\to E$ to the surgery problem with underlying map $f^{-1}(i_K(K))
\to K$ after making $f$ transversal to $i_K \colon K \to E$.

We have shown

\begin{theorem} \label{the:constructing_homotopy_equivalences}
\begin{enumerate}

\item \label{the:constructing_homotopy_equivalences:K}
The diagram~\ref{diagram_for_aut_and_structure_set_for_K_subseteq_E} commutes;

\item \label{the:constructing_homotopy_equivalences:pi}
Each element in the image of $\tau$ is represented by a self-homotopy equivalence
$E \to E$ which induces for some embedded disk $D^{n+k} \subseteq E$ the identity $\id \colon E-D^{n+k} \to E-D^{n+k}$
and in particular induces up to conjugation the identity on the fundamental groups;

\item \label{the:constructing_homotopy_equivalences:S^d}
Suppose additionally that the given map $\phi \colon S^d \to E$ is an embedding.
Let $b' \colon \cals^{\topo}(E) \to \caln(S^d)$ be the map given by making a homotopy equivalence
$f \colon M \to E$ transversal to $\phi$. Then the composite
$$\Omega^{\fr}_{k,k+d}  \xrightarrow{\tau}\cals^{\topo}(E) \xrightarrow{b'} \caln(S^d)$$
is trivial;

\end{enumerate}
\end{theorem}


\typeout{--------------------------------   Section 3 ------------------------------------}

\section{Sphere Bundles}
\label{sec:Sphere_Bundles}

In this section we prove
Theorem~\ref{the:Sphere_bundles_over_surfaces} and
Theorem~\ref{the:Sphere_bundles_over_3_manifolds}.

\begin{proof}
We begin with Theorem~\ref{the:Sphere_bundles_over_surfaces}. It follows from
Theorem~\ref{the:structure_set_for_S^d_bundle_over_K_with_dim(K)le3},
Theorem~\ref{the:constructing_homotopy_equivalences}
\eqref{the:constructing_homotopy_equivalences:K}
and
Theorem~\ref{the:surgery_criterion_for_Borel}
since the $2$-dimensional torus with an appropriate framing yields an element
$\Omega_{2,2+d}^{\fr}$ whose image under the Arf invariant map $\Omega_{2,2+d}^{\fr} \to \bbZ/2$
is non-trivial.
\end{proof}

\begin{proof} Next  we prove Theorem~\ref{the:Sphere_bundles_over_3_manifolds}.
Assertion~\eqref{the:Sphere_bundles_over_3_manifolds:i)} follows from
Theorem~\ref{the:structure_set_for_S^d_bundle_over_K_with_dim(K)le3}
\eqref{the:structure_set_for_S^d_bundle_over_K_with_dim(K)le3:k=3} and
Theorem~\ref{the:surgery_criterion_for_Borel}.

To prove Assertion~\eqref{the:Sphere_bundles_over_3_manifolds:ii)} and
Assertion~\eqref{the:Sphere_bundles_over_3_manifolds:iii)} we use the modified surgery theory from
\cite{Kreck(1999)}, to which we refer for notation. Let $f\colon M \to K \times
S^d$ be a homotopy equivalence. The normal $(d-2)-type$ of $M$ is $K
\times BTOP\langle d-1\rangle$, where $B TOP\langle d-1\rangle$ is the $(d-1)$-connected cover of the
classifying space of topological vector bundles $BTOP$. The reason is that
the restriction of the normal bundle to the $(d-1)$-skeleton, which is
homotopy equivalent to $K$ is determined by $w_2(\nu(M))$. But the
Stiefel-Whitney classes are homotopy invariants, and so $w_2(\nu(M)) = 0$.
If the restriction of the normal bundle to the $(d-1)$-skeleton is
trivial, we obtain a normal $(d-1)$-smoothing of $M$ in $K \times
BTOP\langle d-1\rangle$ by choosing a map from $M$ to $K$ inducing $f_*$ on $\pi_1$ 
up to conjugation and
by choosing a framing on the $(d-1)$-skeleton. Again we use that the
$(d-1)$-skeleton is homotopy equivalent to $K$. Thus a framing on the
restriction of the normal bundle to $K$, considered as the $(d-1)$-skeleton of $M$,
together with $f_* \colon \pi_1(M) \to \pi_1(K)$ determine a normal
$(d-1)$-smoothing of $M$ in $K \times BTOP\langle d-1\rangle $. Since $d \ge 7$,
we conclude that $(d-1)$ is larger than half the dimension of $K \times S^d$.
Thus by \cite[Theorem~3]{Kreck(1999)} and the remark before it and by
\cite[Theorem~4]{Kreck(1999)} and the remark before it the obstruction for replacing a normal
bordism between $K \times S^d$ and $M$ considered as elements of $\Omega
_{d+3}^{TOP\langle d-1\rangle}(K)$ by a $s$-cobordism takes values in the Wall group $L_{d+4}(\pi_1(K))$.
Here we recall that since we assume the Farrell-Jones Conjecture for
$K$-theory we don't need to take the Whitehead torsion into account. From
Theorem~\ref{the:computing_the_structure_set} we know that the $L$-group acts
trivially on the structure set which implies that the obstruction in our
situation vanishes since the action in our situation factors through the
structure set.

Summarizing these considerations we see that $M$ is homeomorphic to $K
\times S^d$ inducing $f_*$ on $\pi_1$ up to conjugation
if and only if  after choosing the
framing on $K$, considered as the $(d-1)$-skeleton of $M$, appropriately the two manifolds are bordant
in $\Omega _{d+3}^{TOP\langle d-1\rangle}(K)$. The different choices of a lift of the normal Gauss
map of $K \times S^d$ to $BTOP\langle d-1\rangle$ correspond to the choice of framings on $K$ and so
$K \times S^d$ is null-bordant for all choices of lifts.
 This implies the following criterion which we will use below:
$M$ is homeomorphic to $K \times S^d$ inducing $f_*$ on $\pi_1$ up to conjugation
if and only if $M$ for one (and
then for all)  topological framings on $K$, considered as the $(d-1)$- skeleton of $M$ , together with the map to
$K$ given by $f_*$ on $\pi_1$ represents the zero class in $\Omega _{d+3}^{TOP\langle d-1\rangle}(K)$.

Next we determine the bordism groups $\Omega _{k}^{TOP\langle d-1\rangle}$ for $d \le k
\le d+3$. If $N$ together with a lift of the normal Gauss map represents
an element in this group we can make it highly connected by surgery. If
$k$ is odd we can even pass to a homotopy sphere which by the topological
Poincar\'e Conjecture is null-bordant. Thus $\Omega _{k}^{TOP\langle d-1\rangle}= 0$ for
$k$ odd (in our range). If $k$ is even, the obstruction for passing to a
homotopy sphere is the signature, if $k = 0 \mod 4$ and the Arf invariant,
if $k = 2 \mod  4$. Since there exist almost parallelizable manifolds with
signature resp. Arf invariant non-trivial, the bordism groups $\Omega
_{k}^{TOP\langle d-1\rangle}$ are $\mathbb Z$, classified by the signature, if $k = 0
\mod 4$ and $\mathbb Z/2$ detected by the Arf invariant, if $k = 2 \mod
4$.

Now we are ready to prove assertion~~\eqref{the:Sphere_bundles_over_3_manifolds:ii)}.
The Atiyah-Hirzebruch spectral sequence
implies for $d= 3 \mod 4$ 
$$\Omega _{d+3}^{TOP\langle d-1\rangle}(K)\cong \Omega
_{d+3}^{TOP\langle d-1\rangle} \oplus H_2(K;\mathbb Z) \cong \mathbb Z/2 \oplus
H_2(K;\mathbb Z).$$ 
Here the first component is determined by the Arf invariant. 
For the detection of second component we note that $H_2(K;\bbZ)$ is isomorphic to $\bbZ$,
and so we can pass to $\mathbb Q$-coefficients. But then the
second component is determined by the higher signatures. Since the Farrell-Jones Conjecture
implies the Novikov Conjecture, the higher signatures of $M \to
K$ agree with the ones of $K \times S^d \to K$ and hence vanish. For the Arf invariant we note that we one can interpret it as
an Arf invariant of a quadratic from given by a Wu-orientation
\cite[Theorem~3.2]{Browder(1969)} and so it is also a homotopy invariant. 
It vanishes for $K \times S^d \to K$ and hence also for $M \to K$.
Thus $M \to K$ is null-bordant in $\Omega _{d+3}^{TOP\langle d-1\rangle}(K)$. 
Hence an application of the criterion above finishes
the proof of
assertion~\eqref{the:Sphere_bundles_over_3_manifolds:ii)}.

To prove assertion~\eqref{the:Sphere_bundles_over_3_manifolds:iii)}
we again use the Atiyah-Hirzebruch spectral sequence to show
$$\Omega _{d+3}^{TOP\langle d-1\rangle}(K)\cong H_1(K;\mathbb Z/2) \oplus
H_3(K;\mathbb Z).$$
We suppose that for $d = 0 \mod 4$ and $d \ge 8$ we have $H_1(K;\mathbb Z/2)\ne
0$. Let $g\colon S^1 \to K$ be a map representing a non-trivial element in
$H_1(K;\mathbb Z/2)$. We consider the connected sum of $K \times S^d$ and
$S^1 \times A$, where $A$ is the framed highly connected topological
manifolds with Arf invariant $1$ (obtained from plumbing two disk bundles
of the tangent bundle of the sphere). So we get normal degree one map
$$\id_{K \times S^d} \# (g \circ p_1) \colon (K \times S^d) \# (S^1 \times A) ~ \to ~ K \times S^d.$$
After composition with the projection $K \times S^d \to K$ we
obtain an element in  $\Omega _{d+3}^{TOP\langle d-1\rangle}(K) = H_1(K;\bbZ/2)
\oplus H_3(K;\bbZ)$. The element is non-trivial since its component
$H_1(K;\bbZ/2)$ is the element represented $g$. This follows from the
product structure of the Atiyah-Hirzebruch spectra sequence.

By the following sequence of surgeries we replace this map by a
homotopy equivalence $f\colon M \to K \times S^d$. Since $\pi_1((K
\times S^d) \# (S^1 \times A)) \cong \pi_1(K \times S^d) \ast
\bbZ$, we can do one $1$-dimensional surgery to make the map an
isomorphism with out changing the homology groups of the universal
coverings up to the middle dimension. Since $A$ is
$(d/2+1)$-connected and $H_{d/2+1}(A)$ is $\bbZ \oplus \bbZ$, the
map is already highly connected and the kernel in the dimension
$(d/2+1)$ on the level of the universal coverings is a free
$\bbZ[\pi_1(K \times S^d)]$-module of rank two. Thus by two
further surgeries we can obtain the desired homotopy equivalence
$f\colon M \to K \times S^d$.
It represents after composition with the projection $K \times S^d \to K$
a non-trivial class in $\Omega _{d+3}^{TOP\langle d-1\rangle}(K)$.
By the criterion mentioned above there is no homeomorphism from $M$ to $K\times S^d$
inducing the same map on $\pi_1$ up to conjugation.
Thus in this situation $K \times S^d$ is not Borel.
\end{proof}


\typeout{--------------------------------   Section 4 ------------------------------------}

\section{Connected Sums}
\label{sec:connected_sums}

In this section we prove Theorem~\ref{the:connected_sums}.
\begin{proof}
The main ingredient is the result of Cappell~\cite[Theorem~0.3]{Cappell(1976a)}
that under our assumptions for every homotopy equivalence $f \colon N \to M_1 \# M_2$
there are $n$-dimensional manifolds
$N_0$ and $N_1$ together with orientation preserving homotopy equivalences
$f_0 \colon N_0 \to M_0$ and $f_1 \colon N_1 \to M_1$ and an orientation preserving homeomorphism
$h \colon N_0 \# N_1 \to N$ such that $f \circ h$ is homotopic to $f_0 \# f_1$.
Now the claim follows from Theorem~\ref{the:surgery_criterion_for_Borel}.
\end{proof}


\typeout{--------------------------------   Section 5 ------------------------------------}

\section{Dimension $3$}
\label{sec:Dimension_3}

Next we prove Theorem~\ref{the:Thurston_and_property_Borel}.
\begin{proof}
If $M$ and $N$ are prime
Haken $3$-manifolds, then every homotopy equivalence $\pi_1(M) \to
\pi_1(N)$ is homotopic to a homeomorphism. This is a result of
Waldhausen (see for instance \cite[Lemma~10.1 and
Corollary~13.7]{Hempel(1976)}). Turaev~\cite{Turaev(1988)} has
extended this result  to showing that a simple homotopy
equivalence between $3$-manifolds with torsionfree fundamental
group is homotopic to a homeomorphism provided that Thurston's
Geometrization Conjecture for irreducible $3$-manifolds with
infinite fundamental group and the $3$-dimensional Poin\-ca\-r\'e
Conjecture are true. This statement remains true if one replaces
simple homotopy equivalence by homotopy equivalence. This follows
from the fact explained below that the Whitehead group of the
fundamental group of a $3$-manifold vanishes provided that
Thurston's Geometrization Conjecture for irreducible $3$-manifolds
with infinite fundamental group is true.

The vanishing of the Whitehead group is proved for Haken manifolds
in Waldhausen~\cite[Theorem~19.3]{Waldhausen(1978b)}. In order to
prove it for prime $3$-manifolds it remains to treat closed
hyperbolic manifolds and closed Seifert manifolds. These cases are
consequences of \cite[Theorem~2.1 on page 263 and
Proposition~2.3]{Farrell-Jones(1993a)}. Now apply the fact that
the Whitehead group of a free amalgamated product is the direct
sum of the Whitehead groups. 
\end{proof}

Every 3-manifold is a generalized topological space form by the following argument.
Suppose that $\pi_1(M)$ is finite. Then the universal covering is a closed
simply connected $3$-manifold and hence homotopy equivalent to $S^3$.
Suppose that $\pi_1(M)$ is infinite. Then the universal covering is a non-compact $3$-manifold and hence homotopy
equivalent to a $2$-dimensional $CW$-complex. This implies that the second homology group of $\widetilde{M}$ 
is a subgroup of a free
abelian group, namely the second chain module of the cellular chain complex of $\widetilde{M}$, and hence free as abelian group
and that all other homology groups of $\widetilde{M}$ are trivial.


\typeout{--------------------------------   Section 6 ------------------------------------}

\section{Dimension $4$}
\label{sec:Dimension_4}

Here we prove Theorem~\ref{the:dimension_4}.
\begin{proof}
\eqref{the:dimension_4:finite_cyclic} Hambleton-Kreck~\cite[Theorem C]{Hambleton-Kreck(1993c)}
show that the homeomorphism type
(including a reference map $M \to B\pi_1(M)$) is determined for a $4$-manifold
with Spin structure and finite cyclic fundamental group by the intersection form
on $M$. Hence such a manifold is Borel.

Here we use the result taken from \cite[10.2B]{Freedman-Quinn(1990)}
that for a $4$-manifold $M$ with Spin structure its signature is divisible by $16$
and its Kirby Siebenmann invariant can be read off from the signature by
$\ks (M) = \sign(M)/8 \mod 2$ and hence
is an invariant of its oriented homotopy type.

Now suppose that $M$ is simply connected and admits no Spin structure. Then there exists
another simply connected $4$-manifold $\ast M$ with the same intersection form but different Kirby Siebenmann invariant
(see~\cite[10.1]{Freedman-Quinn(1990)}.
In particular  $M$ and $\ast M$ are not homeomorphic but they are oriented homotopy equivalent by
\cite{Milnor(1958b)}.
\\[1mm]
\eqref{the:dimension_4:flat}
The Borel Conjecture is true for a flat smooth Riemannian $4$-manifold
(see \cite[page 263]{Farrell(2002)}. Hence such a manifold is strongly Borel.
The manifold $S^1 \times S^3$ is strongly Borel by 
Theorem~\ref{the:S^k_times_S^d}~\eqref{the:S^k_times_S^d:strong_Borel}.

The claim about the connected sum $M \# N$ follows from the version
of Theorem~\ref{the:connected_sums} for dimension $4$ whose proof goes through in dimension
$4$ since the fundamental group of $M \# N$ is good in the sense of Freedman
~\cite{Freedman(1983)}.\\[1mm]

\end{proof}


\typeout{--------------------------------   Section 7 ------------------------------------}

\section{Products of Two Spheres}
\label{sec:Products_of_Two_Spheres}

In this section we give the proof of Theorem~\ref{the:S^k_times_S^d}

\begin{proof}
\eqref{the:S^k_times_S^d:strong_Borel}
If $k + d = 2$, then the property strongly Borel follows
from classical results.

If $k = 1$ and $d \ge 3$ or $k \ge 3$ and $d = 1$, the claim follows from
Theorem~\ref{the:structure_set_for_S^d_bundle_over_K_with_dim(K)le3}
and Theorem~\ref{the:surgery_criterion_for_Borel} since the Farrell-Jones Conjecture is known to be true for
$\pi = \bbZ$.

It remains to treat the case $k,d \ge 2$ and $k +d \ge 4$. Then $S^k \times S^d$ is simply-connected.
The structure set  can be computed by
\begin{eqnarray}
a_1 \times a_2 \colon \cals^{\topo}(S^k \times S^d) & \xrightarrow{\cong} &
 L_k(\bbZ) \oplus L_d(\bbZ),
\label{structure_set_of_S^k_times_S^d}
\end{eqnarray}
where $a_1$ and $a_2$ respectively send the class of an
orientation preserving homotopy equivalence $f \colon M \to S^k
\times S^d$ to the surgery obstruction of the surgery problem with
target $S^k$ and $S^d$ respectively which is obtained from $f$ by
making it transversal to $S^k \times \pt$ and $\pt \times S^d$
respectively. The proof of this claim is analogous to the one of
Theorem~\ref{the:structure_set_of_S^d-bundles_over_K}
(see also \cite[Example~20.4  on page~211]{Ranicki(1992)}).
Hence the structure set is trivial if and only if
both $k$ and $d$ are odd. Now the claim follows from
Theorem~\ref{the:surgery_criterion_for_Borel}. Notice that we can apply surgery theory also
to the $4$-dimensional manifold $S^2 \times S^2$ since its fundamental group is good in the sense of Freedman
~\cite{Freedman(1983)}.
\\[2mm]
\eqref{the:S^k_times_S^d:S^1timesS^2}
This follows from
Theorem~\ref{the:structure_set_for_S^d_bundle_over_K_with_dim(K)le3}~\ref{the:structure_set_for_S^d_bundle_over_K_with_dim(K)le3:k=1}
and Theorem~\ref{the:surgery_criterion_for_Borel} in the case $k \ge 3$. It remains to treat the case
$k = 2$.

Suppose that $S^1 \times S^2$ is Borel. Let $N$ be a homotopy $3$-sphere.
There exists an orientation preserving homotopy equivalence $f \colon S^1 \times S^2 \# N \to S^1 \times S^2$.
Since $S^1 \times S^2$ is Borel by assumption, we can choose $f$ to be a homeomorphism.
By the uniqueness of the prime decomposition $N$  must be homeomorphic to $S^3$.
Hence the $3$-dimensional Poincar\'e Conjecture is true.

Now suppose that the $3$-dimensional Poincar\'e Conjecture is true.
By the Sphere Theorem \cite[Theorem 4.3]{Hempel(1976)}, an irreducible
(closed) $3$-manifold is aspherical if and only if it has infinite fundamental group.
A prime $3$-manifold is either irreducible or is homeomorphic to
$S^1 \times S^2$ \cite[Lemma 3.13]{Hempel(1976)}. Hence for a prime
$3$-manifold $M$  with infinite fundamental group the following statements are equivalent:
i.) M is irreducible, ii.) $M$ is aspheric<al, iii.) $\pi_1(M)$ is not isomorphic to $\bbZ$
and iv.) $M$ is not homeomorphic to $S^1 \times S^2$. Now the prime decomposition of $3$-manifolds implies
that any $3$-manifold with fundamental group $\bbZ$ is homeomorphic to $S^1 \times S^2$.
Hence it suffices to show that any orientation preserving
homotopy equivalence $f \colon S^1 \times S^2 \to S^1 \times S^2$ is homotopic to a homeomorphism.

Since there exists orientation reversing homeomorphisms $S^1 \to S^1$ and
$S^2 \to S^2$, it suffices to treat the case, where $f$ induces the identity on $\pi_1(S^1 \times S^2)$.
Then one can change $f$ up to homotopy that $f$ compatible with the projection
$S^1 \times S^2$, in other words, $f$ is a fiber homotopy equivalence of
the trivial bundle $S^1 \times S^2 \to S^1$. It remains to show
that it is fiber homotopy equivalent to an isomorphism of $S^2$-bundles with structure group
$SO(3)$ over $S^1$. This boils down to showing that the obvious map
$SO(3) \to SG(2)$ is $1$-connected. Analogously to the argument appearing in the proof
of Theorem~\ref{the:motivating_problem:homotopy_type_of_widetilde{M}_for_n_le_5}, but one dimension lower,
one shows that it suffices to show that the (unstable) $J$-homomorphism
$J' \colon \pi_1(SO(2)) \to \pi_3(S^2)$ is bijective. By the Pontrjagin-Thom construction
we obtain a bijection $\pi_3(S^2) \xrightarrow{\cong} \Omega_{1,3}^{\fr}$.
Its composition with $J'$ sends an element in $\pi_1(SO(2))$ to $S^1 \subseteq \bbR^3$ with the
induced framing and hence is surjective. Since $J'$ is a surjective homomorphism of infinite cyclic groups, it is bijective.
\\[1mm]
\eqref{the:S^k_times_S^d:S^2timesS^2}
We have already shown in~\eqref{the:S^k_times_S^d:strong_Borel}  that
$S^2 \times S^2$ is not strongly Borel. It is Borel by
Theorem~\ref{the:dimension_4}~\eqref{the:dimension_4:finite_cyclic}.
\\[1mm]
\eqref{the:S^k_times_S^d:necessary_and_sufficient_for_Borel}
Suppose that $M$ is a Borel manifold. We have to check that conditions
\eqref{the:S^k_times_S^d:necessary_and_sufficient_for_Borel:values_for_k_and_d},
\eqref{the:S^k_times_S^d:necessary_for_Borel:Arf_invariant_for_k}
and~\eqref{the:S^k_times_S^d:necessary_for_Borel:Arf_invariant_for_d} hold.

The $L$-class of $S^k \times S^d$ is concentrated in dimension $0$. We conclude from
Theorem~\ref{the:A_necessary_condition_for_being_a_Borel_manifold} that $S^k \times S^d$ can only be a Borel manifold
if $H_{k+d-4i}(S^k \times S^d;\bbQ)$ is trivial for $i \in \bbZ, i \ge 1$. This implies that neither $k$ nor $d$ are divisible by
four, i.e. condition~\eqref{the:S^k_times_S^d:necessary_and_sufficient_for_Borel:values_for_k_and_d} holds.

Suppose that $k$ is even. By~\eqref{structure_set_of_S^k_times_S^d}
and Theorem~\ref{the:surgery_criterion_for_Borel} there exists an orientation preserving
self-homotopy equivalence $f \colon S^k \times S^d \to S^k \times S^d$ which is transverse to $S^k \times \pt$
such that the Arf invariant of the induced surgery problem $f^{-1}(S^k \times \pt) \to S^k \times \pt$ is
non-trivial. We claim that there exists an orientation preserving homeomorphism
$h \colon S^k \times S^d \to S^k \times S^d$ such that $f$ and $h$ induces the same isomorphism
on $H_n(S^k \times S^d)$ for all $n \in \bbZ$. This is obvious by the K\"unneth formula in the case $k \not = d$ since there
exists a homeomorphism $S^n \to S^n$ of degree $-1$ for all $n \ge 1$. In the case $k = d$, one has to take into
account that $H_k(f)$ respects the intersection form on $H_k(S^k \times S^k) = H_k(S^k) \oplus H_k(S^k) = \bbZ \oplus \bbZ$
which is given by $(x_1,x_2) \cdot (y_1,y_2) \mapsto x_1y_2 + x_2y_1$. This implies that $H_k(f)$ is given by one of the following
matrices
$$
\left(\begin{array}{cc} 1 & 0 \\ 0 & 1 \end{array}\right)
\quad
\left(\begin{array}{cc} -1 & 0 \\ 0 & -1 \end{array}\right)
\quad
\left(\begin{array}{cc} 0 & -1 \\ 1 & 0 \end{array}\right)
\quad
\left(\begin{array}{cc} 0 & 1 \\ -1 & 0 \end{array}\right).
$$
Hence we can find the desired $h$. Now let $g_k \colon S^k \times S^d \to S^k$ be the composition
$\pr_{S^k} \circ f \circ h^{-1}$ for the projection $\pr_{S^k} \colon S^k \times S^d \to S^k$.
Obviously the Arf invariant of the codimension $k$ surgery problem obtained
by making $g_k$ transversal to $\pt \subseteq S^k$ is one and
$g_k$ restricted to $S^k \times \pt$ defines an orientation preserving homotopy equivalence
$S^k \times \pt \to S^k$.
The proof that~\eqref{the:S^k_times_S^d:necessary_for_Borel:Arf_invariant_for_d} holds is
completely analogous.

Now suppose that conditions
\eqref{the:S^k_times_S^d:necessary_and_sufficient_for_Borel:values_for_k_and_d},
\eqref{the:S^k_times_S^d:necessary_for_Borel:Arf_invariant_for_k}
and~\eqref{the:S^k_times_S^d:necessary_for_Borel:Arf_invariant_for_d}
are satisfied. We have to show that $S^k \times S^d$ is Borel what we will do by verifying the
criterion appearing in~Theorem~\ref{the:surgery_criterion_for_Borel}.
In view of assertion~\eqref{the:S^k_times_S^d:strong_Borel}
and condition~\eqref{the:S^k_times_S^d:necessary_and_sufficient_for_Borel:values_for_k_and_d}
we only have to deal with the case, where $k = 2 \mod 4$ or $d = 2 \mod 4$.
We will only explain the most difficult case, where
both $k = 2 \mod 4$ and  $d = 2 \mod 4$ hold, the easier cases, where
$k = 2 \mod 4$ and $d$ is odd or where $d = 2 \mod 4$ and $k$ is odd are then obvious.

Let $g_k \colon S^k \times S^d \to S^k$ be the map appearing in condition
\eqref{the:S^k_times_S^d:necessary_for_Borel:Arf_invariant_for_k}. Define
$f_1 \colon S^k \times S^d \to S^k \times S^d$ to be $g_k \times \pr_{S^d}$
for $\pr_d \colon S^k \times S^d \to S^d$ the projection. This is an orientation preserving
selfhomotopy equivalence satisfying $a_1([f_1]) = 1$ and
$a_2([f_1]) = 0$. Using condition\eqref{the:S^k_times_S^d:necessary_for_Borel:Arf_invariant_for_d}
we construct an orientation preserving
selfhomotopy equivalence $f_2 \colon S^k \times S^d \to S^k \times S^d$ satisfying $a_1([f_2]) = 0$ and
$a_2([f_2]) = 1$.

Obviously we can arrange that $f_1$ and $f_2$  induce the identity on $S^k \vee S^d$.
This together with the identification
$$\sigma_{k+d} \colon \cals^{\topo}(S^k \times S^d) \cong \ker\left(H_{k+d}(S^k \times S^d,\bfL\langle 1 \rangle) \to
H_{k+d}(\pt,\bfL\langle 1 \rangle)\right)$$
implies that the induced map on $\cals^{\topo}(S^k \times S^d)$ the identity. Hence $f_3 = f_1 \circ f_2$
is the desired map because of the formula which has been communicated to us by Andrew Ranicki~\cite{Ranicki(2005z)}
$$[f_2 \circ f_1] ~ = ~ (f_2)_*([f_1]) + [f_2] ~ = ~ [f_1]) + [f_2]$$
and the fact the isomorphism
$$a_1 \times a_2 \colon \cals^{\topo}(S^k \times S^d) \xrightarrow{\cong}  L_k(\bbZ) \oplus L_d(\bbZ),$$
is compatible with the abelian group structures. This finishes the proof of
Theorem~\ref{the:S^k_times_S^d}.
\end{proof}


\typeout{--------------------------------   Section 8 ------------------------------------}

\section{On the Homotopy Type of Certain Low-Di\-men\-sio\-nal Manifolds}
\label{sec:On_the_Homotopy_Type_of_Certain_Low-Dimensional_Manifolds}

We first compute the homology of  the universal covering $\widetilde{M}$
for a manifold appearing in Problem~\ref{pro:Classification_of_certain_low-dimensional_manifolds}.

\begin{lemma} \label{lem:motivating_problem:homology_of_widetilde{M}}
Let $M$ and $K$ be manifolds as described in
Problem~\ref{pro:Classification_of_certain_low-dimensional_manifolds}.
Then $n \ge 4$, if $k = 1$ and $n \ge 5$, if $ k = 2$. Moreover
\begin{eqnarray*}
H_p(\widetilde{M};\bbZ) & \cong_{\bbZ \pi} &
\left\{
\begin{array}{lll}
\bbZ & & p = 0, n-k;
\\
0 & & 1 \le p \le 2;
\\
0 & & p \ge n-2, p \not = n-k,
\end{array}
\right.
\end{eqnarray*}
where $\bbZ$ carries the trivial $\pi$-action. The $\bbZ\pi$-module
$H_3(\widetilde{M})$ is finitely generated stably free if $n = 6$.
\end{lemma}
\begin{proof}
Let $f \colon M \to B\pi$ be the classifying map for $\pi = \pi_1(M)$.
In the sequel we identify $\pi = \pi_1(M) = \pi_1(K)$.
The map is $3$-connected because of $\pi_2(M) = 0$. Let
$\widetilde{f} \colon \widetilde{M} \to \widetilde{K}$ be the induced $\pi$-equivariant
map on the universal coverings. The induced $\bbZ\pi$-chain map
$$C_*(\widetilde{f}) \colon C_*(\widetilde{M}) \to C_*(\widetilde{K})$$
is homological $3$-connected by the Hurewicz Theorem. This implies that its mapping cone
is chain homotopy equivalent to a $\bbZ\pi$-chain complex whose chain modules are
trivial in dimensions $\le 3$. Therefore we obtain isomorphisms
$$H_p(f;\bbZ) \colon H_p(M;\bbZ) ~ \xrightarrow{\cong} H_p(B\pi;\bbZ)$$
and
$$H_p(\widetilde{f};\bbZ) \colon H_p(\widetilde{M};\bbZ) ~ \xrightarrow{\cong} H_p(E\pi;\bbZ)$$
for $p \le 2$ and the induced $\bbZ\pi$-chain map
$$C^{n-*}(\widetilde{f}) \colon C^{n-*}(\widetilde{K}) \to C^{n-*}(\widetilde{M})$$
induces isomorphism for $n-p \le 2$
$$H_p(C^{n-*}(\widetilde{f})) \colon H_p(C^{n-*}(\widetilde{K})) ~ \xrightarrow{\cong}
H_p(C^{n-*}(\widetilde{M})).$$

Obviously $H_p(\widetilde{K};\bbZ) =
H_p(C_*(\widetilde{K}))$ is $\bbZ\pi$-isomorphic to the trivial $\bbZ \pi$-module for $p =
0$ and is trivial for $p \not= 0$. Recall that we have the Poincar\'e $\bbZ\pi$-chain homotopy equivalences
\begin{eqnarray*}
- \cap [M] \colon C^{n-*}(\widetilde{M}) & \to & C_*(\widetilde{M});
\\
- \cap [K] \colon C^{k-*}(\widetilde{K}) & \to & C_*(\widetilde{K}).
\end{eqnarray*}
By Poincar\'e duality applied to $K$ we conclude that
$H_p(C^{n-*}(\widetilde{K}))$ is $\bbZ\pi$-iso\-mor\-phic to $\bbZ$ for
$n-p = k$ and is trivial for $n-p \not= k$.  Hence
$H_p(C^{n-*}(\widetilde{M}))$ is $\bbZ\pi$-isomorphic to $\bbZ$ for
$n-p = k$ and is trivial for $n-p \in \{0,1,2\}, n-p \not = k$.  From
Poincar\'e duality applied to $\widetilde{M}$ we conclude that
$H_p(\widetilde{M};\bbZ)$ is $\bbZ\pi$-isomorphic to $\bbZ$ for $p =
n-k$ and is trivial for $p\in \{n,n-1,n-2\}, p \not = n-k$. We already
know that $H_p(\widetilde{M};\bbZ) = H_p(\widetilde{K};\bbZ) = 0$ for
$p = 1,2$. This implies $n-k \not\in \{1,2\}$ since $\bbZ$ is not
trivial. Hence $n \ge 4$ if $k = 1$ and $n \ge 5$ if $k = 2$. It
remains to show that the $\bbZ\pi$-module $H_3(\widetilde{M})$ is free
if $n = 6$. Let $\cone(f_*)$ be the mapping cone of
$C_*(\widetilde{f}) \colon C_*(\widetilde{M}) \to C_*(\widetilde{K})$.
The composition
$$C_*(\widetilde{f}) ~ \circ ~ \left(- \cap [M]\right) ~ \circ ~ C^{6-*}(\widetilde{f})
\colon C^{6-*}(\widetilde{K}) \to C_*(\widetilde{K})$$
is obviously zero. Hence the $\bbZ\pi$-chain map
$$\left(- \cap [M]\right) ~ \circ ~ C^{6-*}(\widetilde{f})
\colon C^{6-*}(\widetilde{K}) \to C_*(\widetilde{M})$$
induces in the obvious way a $\bbZ\pi$-chain map
$$C^{6-*}(\widetilde{K}) \to \Sigma^{-1} \cone(f_*).$$
Let $D_*$ be its mapping cone.
It inherits from the structure of a symmetric $6$-dimensional Poincar\'e chain complex on $C_*(\widetilde{M})$
the  structure of a symmetric Poincar\'e chain complex on $D_*$ with Poincar\'e dimension $6$.
This follows from \cite[Proposition~4.1 on page~141]{Ranicki(1973a)}.
The homology of $D_*$ is zero in dimensions different from $3$ and $H_3(D_*)$ is
$\bbZ\pi$-isomorphic to $H_3(\widetilde{M})$. Let $E_*$ be the $\bbZ\pi$-subchain complex of $D_*$
for which $E_k = 0$ for $k \le 2$, $E_3 = \ker(d_3 \colon D_3 \to D_2)$ and $E_k = D_k$ for
$k \ge 4$. Since
$$0 \to \ker(d_3) \to D_3 \xrightarrow{d_3} D_2 \xrightarrow{d_2} D_1  \xrightarrow{d_1} D_0 \to D_{-1} \to 0$$
is exact, $E_*$ is a finitely generated free  $\bbZ\pi$-chain
complex. The inclusion $i \colon E_* \to D_*$ is a homology
equivalence and hence a $\bbZ\pi$-chain homotopy equivalence. In
particular we can pullback the  structure of a symmetric
Poincar\'e chain complex of Poincar\'e dimension $6$ on $D_*$ to
$E_*$. Hence $E^{6-*}$ is a $\bbZ\pi$-chain complex which is
concentrated in dimensions $-1$, $0$, $1$, $2$ and $3$, whose
homology is zero in dimensions different from $3$ and for which
$H_3(E^{6-*})$ is $\bbZ\pi$-isomorphic to $H_3(\widetilde{M})$.
This implies that there is an exact sequence of $\bbZ\pi$-modules
$$0 \to H_3(\widetilde{M})  \to E_3^* \to E_4^* \to E_5^* \to E_6^* \to E_7^* \to 0.$$
Hence $H_3(\widetilde{M})$ is a finitely generated
stably free $\bbZ\pi$-module.
This finishes the proof
of Lemma~\ref{lem:motivating_problem:homology_of_widetilde{M}}.
\end{proof}

Next we determine the homotopy type of $M$ in the case $n = 5$.

\begin{theorem} \label{the:motivating_problem:homotopy_type_of_widetilde{M}_for_n_le_5}
Let $M$ and $K$ be as in
Problem~\ref{pro:Classification_of_certain_low-dimensional_manifolds}.
Let $f \colon M \to K = B\pi$ be the classifying map
for the $\pi$-covering $\widetilde{M} \to M$.
Suppose $n \le 5$. Then

\begin{enumerate}

\item
Let $k = 2$. The homotopy type of $M$ is determined by its second Stiefel-Whitney class.
Namely, there is precisely one fiber bundle $S^3 \to E \xrightarrow{p} K$ with structure group $SO(4)$
whose second Stiefel-Whitney class agrees
under the isomorphism $H^2(f;\bbZ/2) \colon H^2(K;\bbZ/2) \xrightarrow{\cong}
H^2(M;\bbZ/2)$ with the second Stiefel-Whitney class of $M$. There exists a homotopy equivalence $g \colon M \to E$
such that $p \circ g$ and $f$ are homotopic;

\item Suppose that $k=2$ and the second Stiefel-Whitney class of $M$ vanishes, or that $ k = 1$.
Then there is a homotopy equivalence $g \colon M \to S^{n-k} \times B\pi$
such that $\pr_{B\pi} \circ g$ and $f$ are homotopic.
\end{enumerate}
\end{theorem}
\begin{proof}
Obviously the second assertion is a special case of the first one.
The first one is proven as follows.

We conclude from Lemma~\ref{lem:motivating_problem:homology_of_widetilde{M}} that
for $n \le 5$ the universal covering $\widetilde{M}$ has the same homology as $S^{n-k}$ and that $\pi$
acts trivially on the homology of $\widetilde{M}$. Recall that $K$ is a model for $B\pi$.
There is a fibration $F \to E \xrightarrow{p} K$ together with a homotopy
equivalence $g \colon E \to M$ such that $f \circ g \simeq p$ and $F$ is homotopy
equivalent to $\widetilde{M}$. Hence the fibration $p \colon E \to K$ has $S^{n-k}$ as
fiber.  Since the $\pi$-action on the homology of $\widetilde{M}$ is trivial,
this spherical fibration is orientable, i.e.
the fiber transport of this fibration $\pi = \pi_1(K) \to [F,F]$ is trivial.

Such fibrations over $K$ are classified by maps $u \colon K \to BSG(n-k)$,
where $BSG(n-k)$ is the classifying space of the monoid $SG(n-k)$ of
orientation preserving self-homotopy equivalences $u \colon S^{n-k} \to S^{n-k}$.
We have $n-k \ge 3$. The space $SG(n-k)$ is connected and hence the space
$BSG(n-k)$ is simply connected.

Suppose that $k = 1$. Each map $S^1 \to BSG(n-k)$ is trivial up to homotopy
and the claim follows.

Now suppose $k = 2$. Then $n \ge 5$ by Lemma~\ref{lem:motivating_problem:homology_of_widetilde{M}}
Hence $n = 5$ because of the assumption $n \le 5$.
There is an obvious commutative diagram
$$
\begin{CD}
SO(3) @>>>  SO(4) @>>> S^3
\\
@VJ' VV @V J  VV @V\id VV
\\
\Omega^3S^3 @>>> SG(3) @>>> S^3
\end{CD}
$$
where the horizontal maps are fibrations and $J$ comes from the obvious
action of $SO(4)$ on $S^3$. The following diagram commutes
\comsquare{\pi_1(SO(3))}{\cong}{\pi_1(SO)}{\pi_1(J')}{J_1}{\pi_1(\Omega^3S^3) = \pi_4(S^3)}{\cong}{\pi_1^s}
where the horizontal maps are given by stabilization and are isomorphisms and
$J_1$ is the $J$-homomorphism. The abelian groups $\pi_1(SO)$ and $\pi_1^s$ are both isomorphic to $\bbZ/2$.
The map $J_1$ is bijective  by \cite{Adams(1966)}.
Hence $J \colon SO(4) \to SG(3)$ is a map of connected space which induced an isomorphism
on the fundamental groups. This implies that
$BJ \colon BSO(4) \to BSG(3)$ is a map of simply connected spaces inducing an isomorphism
on $\pi_2$. Let $w_2 \colon BSG(3) \to K(\bbZ/2,2)$ be given by the second Stiefel-Whitney class.
It and its composite $w_2  \circ BJ \colon BSO(4) \to K(\bbZ/2,2)$ are $3$-connected. Since $K$ is
$2$-dimensional, we conclude that every orientable fibration $S^3 \to E \to K$ is fiber homotopy
equivalent to a fiber bundle $S^3 \to E \to K$ with structure group $SO(4)$ and
two such fiber bundles with structure group $SO(4)$ over $K$ are isomorphic if and  only if
their second Stiefel-Whitney classes agree.

We have $H^2(f;\bbZ/2)(w_2(p)) = w_2(M)$ since for
a vector bundle $\xi \colon E' \to K$ we get a decomposition
$TE'|K \cong \xi \oplus TK$ and the Stiefel-Whitney classes of $K$ are trivial.

This finishes the proof of Theorem~\ref{the:motivating_problem:homotopy_type_of_widetilde{M}_for_n_le_5}.
\end{proof}

In dimension $n = 6$ we can at least compute the homology of the universal
covering up to stable isomorphism.

\begin{lemma} \label{lem:motivating_problem_homology_type_of_widetilde{M}_for_n=6}
Let $M$ and $K$ be as in
Problem~\ref{pro:Classification_of_certain_low-dimensional_manifolds}.
Suppose $n = 6$.
Let $\chi(M)$ be the Euler characteristic.

Then $-\chi(M) + (1 + (-1)^k) \cdot \chi(K) \ge 0$ and we get
\begin{eqnarray*}
H_p(\widetilde{M};\bbZ) & \cong_{\bbZ \pi} &
\left\{
\begin{array}{lll}
\bbZ & & p = 0, n-k;
\\
0 & & 0 \le p \le 2;
\\
0 & & p \ge n-2, p \not = n-k.\\
\end{array}
\right.
\end{eqnarray*}
and for some $s \ge 0$
$$H_3(\widetilde{M};\bbZ)  \oplus \bbZ\pi^s ~ \cong_{\bbZ \pi} ~
\bbZ\pi^{-\chi(M) + (1 + (-1)^k) \cdot \chi(K)} \oplus \bbZ\pi^s.$$
\end{lemma}
\begin{proof}
Because of Lemma~\ref{lem:motivating_problem:homology_of_widetilde{M}} it suffices to
show for the finitely generated stably free $\bbZ\pi$-module
$H_3(\widetilde{M})$ that the classes of $H_3(\widetilde{M})$ and
$\bbZ\pi^{-\chi(M) + (1 + (-1)^k) \cdot \chi(K) }$ agrees in $K_0(\bbZ \pi)$.

We obtain a finite projective resolution for the trivial
$\bbZ\pi$-module $\bbZ$ by $C_*(\widetilde{K})$ and get in $\widetilde{K}_0(\bbZ\pi)$
$$[\bbZ] ~ = ~ \chi(K) \cdot [\bbZ\pi].$$
Therefore every homology module $H_p(\widetilde{M})$ for $p \not= 3$ possesses a finite
projective $\bbZ\pi$-resolution and hence
defines an element in $K_0(\bbZ\pi)$. This automatically implies
that the same is true for $H_3(\widetilde{M})$. We compute in $K_0(\bbZ\pi)$
\begin{eqnarray*}
 [H_3(\widetilde{M})]
& = &
- \sum_{k = 0}^6 (-1)^k \cdot [H_k(\widetilde{M})] + [\bbZ] + (-1)^{6-k}[\bbZ]
\\
& = &
- \sum_{k = 0}^6 (-1)^k \cdot [C_k(\widetilde{M})] + (1 + (-1)^k) \cdot [\bbZ]
\\
& = & -\chi(M) \cdot [\bbZ\pi] + (1 + (-1)^k) \cdot \chi(K) \cdot [\bbZ\pi]
\\
& = & \left( -\chi(M) + (1 + (-1)^k) \cdot \chi(K) \right)  \cdot [\bbZ\pi].
\end{eqnarray*}
Now Lemma~\ref{lem:motivating_problem_homology_type_of_widetilde{M}_for_n=6} follows.
\end{proof}


\typeout{--------------------------------   Section 9 ------------------------------------}

\section{On the Classification of Certain Low-Di\-men\-sio\-nal Manifolds}
\label{sec:On_the_Classification_of_Certain_Low-D-mensional_Manifolds}

\begin{theorem} \label{the:motivating_problem:topological_structure_set}
Let $M$ and $K$ be as in
Problem~\ref{pro:Classification_of_certain_low-dimensional_manifolds}.
Suppose $k= 2$.

Then $n \ge 5$ and there is
an isomorphism of abelian groups
$$\cals^{\topo}(M) \xrightarrow{\cong} \bbZ/2.$$
If $f \colon M \to M$ is an orientation preserving  self-homotopy equivalence and $g \colon N \to M$
is an
orientation preserving homotopy equivalence of manifolds, then we obtain in
$\cals^{\topo}(M)$
$$[g \circ f] ~ = ~ [g] + [f].$$
\end{theorem}
\begin{proof}
Notice for the sequel
that the Farrell-Jones Conjecture for algebraic $K$- and $L$-theory
holds for $\pi_1(K)$ (see~\cite{Farrell-Jones(1993c)}).

The Atiyah-Hirzebruch spectral sequence
yields two commutative diagrams with exact columns respectively rows
\begin{small}
$$\begin{CD}
0 && 0
\\
@VVV @VVV
\\
H_0(M;L_6(\bbZ)) \oplus H_2(M;L_4(\bbZ))
@> H_0(f;L_6(\bbZ)) \oplus H_2(f;L_4(\bbZ)) > \cong >
H_0(K;L_6(\bbZ)) \oplus H_2(K;L_4(\bbZ))
\\
@VVV @VVV
\\
H_6(M;\bfL\langle 1 \rangle) @>H_6(f;\bfL\langle 1 \rangle) >> 
H_6(K;\bfL\langle 1 \rangle)
\\
@Ve_6VV @Ve_6VV
\\
 H_4(M;L_2(\bbZ)) @>>> 0
\\
@VVV @VVV
\\
0 && 0
\end{CD}
$$
\end{small}
and
$$\begin{CD}
0 @>>>  H_1(M;L_4(\bbZ))
@>>> H_5(M;\bfL\langle 1 \rangle)
@>e_5 >> H_3(M;L_2(\bbZ)) @ >>> 0
\\
& & @V H_1(f;L_4(\bbZ)) V \cong V
@V H_5(f;\bfL\langle 1 \rangle) VV
@V 0 VV
\\
0 @>>>  H_1(K;L_4(\bbZ))
@>>> H_5(K;\bfL\langle 1 \rangle)
@>e_5 >> 0 @ >>> 0
\end{CD}
$$
We conclude that the restriction of
$$e_n \colon H_n(M;\bfL\langle 1 \rangle) \to H_{n-2}(M;L_2(\bbZ))$$
to the kernel of
$$H_n(f;\bfL\langle 1 \rangle)) \colon H_n(M;\bfL\langle 1 \rangle) \to H_n(K;\bfL\langle 1 \rangle)$$
is bijective.
We conclude from Theorem~\ref{the:computing_the_structure_set}
that the following composition is an isomorphism for $n = \dim(M)$.
$$\cals^{\topo}(M) \xrightarrow{\sigma_n}
H_n(M;\bfL\langle 1 \rangle) \xrightarrow{e_n} H_{n-2}(M;L_2(\bbZ)).$$
We have the following composition of isomorphisms
$$H^2(K;L_2(\bbZ)) \xrightarrow{H^2(f;L_2(\bbZ))} H^2(M;L_2(\bbZ))
\xrightarrow{- \cap [M]}  H_{n-2}(M;L_2(\bbZ)).$$
Since $k = 2$ by assumption, we get
$$H^2(K;L_2(\bbZ)) \cong L_2(\bbZ) \cong \bbZ/2.$$

If $f \colon M \to M$ is an orientation preserving  self-homotopy equivalence and $g \colon N \to M$
orientation preserving homotopy equivalence of manifolds, then we obtain in
$\cals^{\topo}(M)$ the formlua which has been communicated to us by Andrew Ranicki~\cite{Ranicki(2005z)}
$$[f \circ g] ~ = ~ f_*[g] + [f]$$
Since any automorphism of the group $\bbZ/2$ is the identity,
this finishes the proof of
Theorem~\ref{the:motivating_problem:topological_structure_set}.
\end{proof}

Now we are ready to prove Theorem~\ref{the:motivating_problem_n_le_5} and
Theorem~\ref{the:motivating_problem_n=6}~\eqref{the:motivating_problem_n=6:Borel}.

\begin{proof} In view of Theorem~\ref{the:motivating_problem:homotopy_type_of_widetilde{M}_for_n_le_5},
Theorem~\ref{the:motivating_problem:topological_structure_set} and
Theorem~\ref{the:surgery_criterion_for_Borel} it remains to show that
there exists a homotopy equivalence $f \colon M \to M$ inducing the up to conjugation the identity
on the fundamental group which does not represent the trivial element $\id \colon M \to M$ in
$\cals^{\topo}(M)$. We want to conclude this from
Theorem~\ref{the:constructing_homotopy_equivalences}~\eqref{the:constructing_homotopy_equivalences:K}
and \eqref{the:constructing_homotopy_equivalences:pi}
As soon as we have provided the necessary data, the claim follows from
Theorem~\ref{the:constructing_homotopy_equivalences}~\eqref{the:constructing_homotopy_equivalences:K}
since the $2$-dimensional torus with an appropriate framing yields an element
$\Omega_{2,2+d}^{\fr}$ for $d \ge 1$ whose image under the Arf invariant map $\Omega_{2,2+d}^{\fr} \to \bbZ/2$
is non-trivial.

Since the dimension of $K$ is less or equal to the connectivity of $S^d$,
we can choose a map $i \colon K \to S^d$ such that $p \circ i$ is homotopic to the identity.
By the homotopy lifting property we can arrange that $p \circ i = \id_K$. In particular we see that
$i \colon K \to M$ is an embedding. It remains to construct the map
$\alpha \colon S^d \to M$.

Let $\nu = \nu(i)$ be the normal bundle of $i \colon K \to M$ and
$D\nu$ be the associated disk bundle. Let $j_1 \colon D^d
\subseteq D\nu \subseteq M$ be the inclusion of the fiber of a
point $x \in K$. Obviously $D^d$ intersects $K$ transversally and
the intersection consists of one point. It suffices to show that
the map $j_0 \colon S^{d-1} = \partial D^d \to M-K$ given by
restricting $j_1$ to $S^{d-1}$ is nullhomotopic  because then we
can extend $j_0$ to a map $j_2 \colon D^d \to M-K$ and can define the
desired map $\phi$ by $j_2 \cup_{j_0} j_1 \colon S^d = D^d
\cup_{S^{d-1}} D^d \to M$.

The map $i \colon K \to M$ and the map $p \colon M \to K$ induce isomorphisms on the fundamental groups.
Fix a universal covering $\widetilde{K} \to K$. Its pullback with $p \colon M \to K$ is the universal
covering $\widetilde{M} \to M$. Obviously we get $\pi$-equivariant maps
$\widetilde{i} \colon \widetilde{K} \to \widetilde{M}$ and $\widetilde{p} \colon
\widetilde{M} \to \widetilde{K}$ covering $i$ and $p$ and satisfying
$\widetilde{p} \circ \widetilde{i} = \id_{\widetilde{K}}$.
Let $\widetilde{j_1} \colon D^d \to \widetilde{M}$ and
$\widetilde{j_0} \colon S^{d-1} \to \widetilde{M}$ be lifts of $j_1$ and $j_0$
with $\widetilde{j_1}|_{S^{d-1}} = \widetilde{j_0}$.
The map $\pi_{d-1}(\widetilde{M}-\widetilde{i}(\widetilde{K})) \to \pi_{d-1}(M-i(K))$
is induced by a covering map and hence an isomorphism because of $d-1 \ge 2$.
It sends the class of $\widetilde{j_0}$ to the class of $j_0$. Hence it suffices to show that
the class of $\widetilde{j_0}$ in $\pi_{d-1}(\widetilde{M}-\widetilde{i}(\widetilde{K}))$ is zero.
Its image under the map $\pi_{d-1}(\widetilde{M}-\widetilde{i}(\widetilde{K})) \to
\pi_{d-1}(\pi_{d-1}(\widetilde{M}))$ is zero, a nullhomotopy for the image is given by $\widetilde{j_1}$.
Since $\widetilde{M}$ is simply connected and its homology is trivial in dimensions
$\le d-2$ by Lemma~\ref{lem:motivating_problem:homology_of_widetilde{M}},
the space $\widetilde{M}$ is $(d-2)$-connected. Since the codimension of $\widetilde{i}(\widetilde{K})$ in $\widetilde{M}$ is $d$,
the inclusion $\widetilde{M}-\widetilde{i}(\widetilde{K}) \to \widetilde{M}$ is $(d-1)$-connected and hence
$\widetilde{M}-\widetilde{i}(\widetilde{K})$ is $(d-2)$-connected. Therefore the
Hurewicz homomorphism
$$\pi_{d-1}(\widetilde{M}-\widetilde{i}(\widetilde{K})) \xrightarrow{\cong}
 H_{d-1}(\widetilde{M}-\widetilde{i}(\widetilde{K}))$$
is an isomorphism. The image of the class of $\widetilde{j_0}$ under the Hurewicz homomorphism is send to zero under the
map $H_{d-1}(\widetilde{M}-\widetilde{i}(\widetilde{K})) \to H_{d-1}(\widetilde{M})$.
Hence it remains to show that this map is injective. By the long exact sequence of the pair
it suffices to prove that
$H_d(\widetilde{M}) \to H_d(\widetilde{M},\widetilde{M}-\widetilde{i}(\widetilde{K}))$
is surjective. By Poincar\'e duality we get a commutative diagram with isomorphisms
as vertical maps
\comsquare{H_d(\widetilde{M})}{}{H_d(\widetilde{M},\widetilde{M}-\widetilde{i}(\widetilde{K}))}
{\cong}{\cong}
{H^k(\widetilde{M})}{}{H^k(\widetilde{K})}
where the lower horizontal arrow is induced by $\widetilde{i}$ and is split surjective because
of $\widetilde{p} \circ \widetilde{i} = \id_{\widetilde{K}}$. This finishes the proofs
of Theorem~\ref{the:motivating_problem_n_le_5} and
Theorem~\ref{the:motivating_problem_n=6}~\eqref{the:motivating_problem_n=6:Borel}.
\end{proof}

Next we give the proof of Theorem~\ref{the:motivating_problem_n=6}~\eqref{the:motivating_problem_n=6:Spin_case}.
\begin{proof}
Let $M$ be a closed topological oriented $6$-manifold with $\pi_1(M) \cong \pi_1(K)$,
where $K$ is a $K (\pi ,1)-$manifold of dimension $\le2$, and $\pi_2(M) = 0$.
We want to prove that the $\bbZ\pi$-isomorphism class of the intersection form
on $H_3(\widetilde{M})$ determines the homeomorphism type.
The normal $2$-type of $M$ in the sense of Kreck~\cite{Kreck(1999)}
is
$$K \times BTopSpin \xrightarrow{\pr_2} BTopSpin \xrightarrow{P} BSTOP,$$
if $M$ is a topological $Spin$-manifold and
$$K \times BTopSpin \xrightarrow{E \times P} B STOP \times BSTOP \xrightarrow{\oplus} BSTOP$$
if $M$ is not a topological $Spin$-manifold, where $E$ is a vector bundle over $K$ with
$w_2 (E) \not= 0$. This normal $2$-type  is determined by $\pi = \pi_1(K)$
and $w_2 = w_2 (M) \in \bbZ/2$ and so we denote it by $B(\pi,w_2)$.  Notice that
$(M; \overline{\nu})$ determines a class in the bordism group $\Omega_6(B(\pi,w_2))$.

Now we want to apply \cite[Corollary 3]{Kreck(2004b)}. It says that if $(M'; \overline{\nu})$
is another normal $2$-smoothing in $B (\pi, w_2)$, then $M$ and $M'$ are homeomorphic,
where the homeomorphism is compatible with the maps on $\pi_1$, if and only if the pairs
determine the same class in $\Omega_6(B(\pi,w_2))$ and the intersection form together with quadratic refinement on
$K(\pi_3(M)) \to \pi_3 (B (\pi, w_2)) = \pi _3 (M)$, which by
Lemma~\ref{lem:motivating_problem:homology_of_widetilde{M}}
is stably free and by Poincar\'e duality unimodular, are isomorphic.

Since $\langle  w_4 (B(\pi, w_2)), \pi_4 (B(\pi, w_2))\rangle \not= 0$
we have to consider the quadratic refinement with values in
$\bbZ \pi/\langle x + \overline{x}, 1\rangle$. This is a quadratic
refinement with respect to a form parameter in the sense of Bak~\cite{Bak(1981)}.
Since $\pi$ has no element of order $2$ and the quadratic form on $\pi_3(M)$ takes
values in $\bbZ [\pi]/\langle a + \overline{a}, 1\rangle = 0$, the quadratic form is determined by the intersection form.
So it remains to show that the intersection form determines the bordism class in
$\Omega_6 (B(\pi,w_2))$.

The Atiyah-Hirzebruch spectral sequence shows
that $\Omega_6 (B(\pi ,w_2)) \cong H_2 (K; \Bbb Z))$.  If $\dim(K) = 1$, then
$\Omega_6 (B(\pi ,w_2))$ is trivial and the claim follows.

Suppose $\dim(K) = 2$. Then $\Omega_6 (B(\pi, w_2))$ is isomorphic to
$\bbZ$.  Consider the following composite
$$\Omega_6 (B(\pi ,w_2)) \xrightarrow{\alpha} \Omega_6(B\pi) \xrightarrow{\sigma} L^6(\bbZ \pi)$$
where $\alpha$ is induced by the obvious map $B(\pi,w_2) \to B\pi$ and the second map
is given by the symmetric signature in the sense of
Ranicki~\cite{Ranicki(1980a)}, \cite{Ranicki(1981)}. The argument in the proof of
Lemma~\ref{lem:motivating_problem:homology_of_widetilde{M}} shows that the class
$\beta \circ \alpha([M; \overline{\nu}]$ which is given by the chain complex of the
associated $\pi$-covering of $M$ and the Poincar\'e chain homotopy equivalence is the same as the class
represented by the intersection form, since the two are obtained from one another by algebraic surgery.
Hence the intersection from determines $\beta \circ \alpha([M; \overline{\nu}]$. Since
$\Omega_6 (B(\pi,w_2))$ is torsionfree, it remains to show that
$\beta \circ \alpha$ is rationally injective.

There is the following commutative diagram
of $\bbZ$-graded rational vector spaces
as explained in~\cite[page~728]{Kreck-Leichtnam-Lueck(2002)}
$$
\begin{CD}
\left(\Omega_*(B\pi)\otimes_{\Omega_*(*)} \bbQ\right)_n
@> \overline{D} > \cong >KO_n(B\pi)\otimes_{\bbZ} \bbQ @>i>>
K_n(B\pi) \otimes_{\bbZ} \bbQ
\\
@V \sigma VV @VA_{\bbR}  VV @V A_{\bbC} VV
\\
L^n(C_r^*(\pi;\bbR))\otimes_{\bbZ} \bbQ @> \sign > \cong >
KO_n(C_r^*(\pi;\bbR))\otimes_{\bbZ} \bbQ @>j>>
K_n(C_r^*(\pi)) \otimes_{\bbZ} \bbQ
\end{CD}
$$
The map $A_{\bbC}$ is the assembly map appearing in the Baum Connes Conjecture and
is known to be an isomorphism. The map $i$ is a change of rings map and known to be
rationally injective.  The map $\sigma$ factorizes as
$$\left(\Omega_*(B\pi)\otimes_{\Omega_*(*)} \bbQ\right)_n \to L^n(\bbZ \pi)\otimes_{\bbZ} \bbQ \to L^n(C_r^*(\pi;\bbR))\otimes_{\bbZ} \bbQ.$$
Hence it suffices to show that the composite
\begin{multline*}
\Omega_6 (B(\pi ,w_2)) \otimes_{\bbZ} \bbQ \xrightarrow{\alpha \otimes_{\bbZ} \id_{\bbQ}} \Omega_6(B\pi) \otimes_{\bbZ} \bbQ
\to \Omega_*(B\pi)\otimes_{\Omega_*(*)} \bbQ
\to KO_6(B\pi)\otimes_{\bbZ} \bbQ
\end{multline*}
is injective. This follows from a spectral sequence argument.
This shows that the intersection form determines the homeomorphism type.

The next question is which unimodular forms $\lambda$ on stably free $\mathbb ZÊ[\pi_1]-$modules $V$
can be realized as intersection forms of manifolds under consideration. Since, if $(V, \lambda)$
can be realized, and $(V, \lambda)$ splits off a hyperbolic form, i.e.
$(V, \lambda) = (V', \lambda') \perp H$, then $(V', \lambda')$ can be realized
by surgery on the hyperbolic plane, the realization problem is reduced to the
stable realization problem: Which elements in $\widetilde L_6 (\pi_1 (K))$ can be realized by a stable homeomorphism class?

If $w_2 = 0$, the answer is: All. The reason is that we have a commutative diagram
$$
\begin{array}{ccccc }
   \Omega^{TopSpin}_4   & \stackrel{\cong}{\longrightarrow} & L_4 (\{e\})    \\
   &&\\
     \cong \; \downarrow \; \times \: K & &   \downarrow \; \cong   \\
     &&\\
     \Omega^{TopSpin}_6 (K) & \longrightarrow &\widetilde  L_6 (\pi_1 (K))
\end{array}
$$
If $w_2 \not= 0$ we don't know the answer.
\end{proof}


\typeout{--------------------------------   Section 10 ------------------------------------}

\section{A Necessary Condition for Being a Borel Manifold}
\label{sec:A_Necessary_Condition_for_Being_a_Borel_Manifold}

Next we prove Theorem~\ref{the:A_necessary_condition_for_being_a_Borel_manifold}.
\begin{proof}
For every homology theory satisfying the disjoint union axiom and
hence in particular for $H_*(-;\bfL\langle 1\rangle)$ there is a
natural Chern character (see Dold~\cite{Dold(1962)}).
\begin{eqnarray}
\ch_n(X) \colon \bigoplus_{i \in \bbZ, i \ge 1} H_{n-4i}(X;\bbQ) \xrightarrow{\cong} H_n(X;\bfL\langle 1\rangle) \otimes_{\bbZ} \bbQ
\label{chern_character_for H_*(-;l<1>)}
\end{eqnarray}
For an $n$-dimensional manifold $M$ the composite
$$\cals^{\topo}(M) \xrightarrow{\sigma_n} H_n(M;\bfL\langle 1\rangle) \to
H_n(M;\bfL\langle 1\rangle) \otimes_{\bbZ} \bbQ \xrightarrow{\ch_n(M)^{-1}}
\bigoplus_{i \in \bbZ, i \ge 1} H_{n-4i}(M;\bbQ) $$
sends the class of an orientation preserving homotopy equivalence
$f \colon N \to M$ to the element $\{f_*(\call(N) \cap [N]) - \call(M) \cap [M]\mid i \in \bbZ, i \ge 1\}$.
(see \cite[Example~18.4 on page~198]{Ranicki(1992)}).

Now suppose that $M$ is a Borel manifold. Then by Theorem~\ref{the:surgery_criterion_for_Borel}
the operation of the group of homotopy classes
of self-homotopy equivalence $M \to M$ which induce the identity on the fundamental group up to conjugation
on the topological structure set $\cals^{\topo}(M)$ is transitive.  This implies
that this group acts also transitive on the image of
the composite
$$\cals^{\topo}(M) \xrightarrow{\sigma_n} H_n(M;\bfL\langle 1\rangle)  \to
H_n(M;\bfL\langle 1\rangle) \otimes_{\bbZ} \bbQ \xrightarrow{\ch_n(M)^{-1}}
\bigoplus_{i \in \bbZ, i \ge 1} H_{n-4i}(M;\bbQ).$$
This image is obviously an abelian subgroup of $\bigoplus_{i \in \bbZ, i \ge 1} H_{n-4i}(N;\bbQ)$
and agrees with the set
$$S: = \{f_*(\call(M) \cap [M]) - \call(M) \cap [M] \mid [f] \in \hoaut_{\pi}(M)\}.$$
By the exactness of the surgery sequence the $\bbQ$-submodule
generated by $S$ contains the kernel of the map induced by the classifying map $c \colon M \to B\pi$
$$c_* \colon \bigoplus_{i \in \bbZ, i \ge 1} H_{n-4i}(M;\bbQ)  \to
\bigoplus_{i \in \bbZ, i \ge 1} H_{n-4i}(B\pi;\bbQ) .$$
\end{proof}

Now we are ready to prove Theorem~\ref{the:necessary_condition_for_sphere bundles}.
\begin{proof}
We first show that the set of homotopy classes of orientation preserving homotopy equivalences $f \colon E \to E$ which induce
up to conjugation the identity on the fundamental group is finite. Since $K$ is aspherical and
$p \colon E \to K$ induces an isomorphism on the fundamental groups because of $d \ge 2$,
the maps $p$ and $p \circ f$ are homotopic. By the homotopy lifting property we can assume
that $p \circ f = p$ holds. Hence it suffices to show that the set of fiber homotopy
classes of fiber homotopy equivalences $f \colon E \to E$ which cover the identity $\id_K \colon K \to K$
and induce a map of degree one on the fibers is finite. By elementary obstruction theory this follows
if the $i$-th homotopy group of the space $SG(S^d)$ of self-maps $S^d \to S^d$ of degree one is finite for $i \le k$.
There is an obvious fibration $\Omega^dS^d \to SG(S^d) \to S^d$. The long exact homotopy sequence yields the exact sequence
$\pi_{i+d}(S^d) \to \pi_i(SG(S^d)) \to \pi_i(S^d)$, where we take the obvious base points.
If $d$ is odd, $\pi_j(S^d)$ is finite for all $j \ge 0$, and, if $d$ is even, $\pi_j(S^d)$ is finite for $j \le 2d-2$. This
has been proved by Serre~\cite{Serre(1953)}. Hence  $\pi_i(SG(S^d))$ is finite if $i \ge 0$ and
$d$ is odd or if $i \le d-2$.

Since $K$ is aspherical and there is a map $i \colon K \to E$ with $p \circ i \simeq \id_K$,
the kernel of the map $c_* \colon H_{k+d-4i}(E) \to H_{k+d-4i}(B\pi_1(E);\bbQ)$ induced by the classifying map
$c = p \colon E \to B\pi_1(E) = K$ contains  $H_{k-4i}(K;\bbQ)$. Now the claim follows from
Theorem~\ref{the:A_necessary_condition_for_being_a_Borel_manifold} because the abelian subgroup
$S$ of a $\bbQ$-module appearing there is finite and hence trivial.
\end{proof}


\typeout{--------------------------------   Section 11 ------------------------------------}

\section{Integral Homology Spheres}
\label{sec:Integral_Homology_Spheres}

In this section we prove
Theorem~\ref{the:Integral_homology_spheres}
\begin{proof}
\eqref{the:Integral_homology_spheres:L-criterion}
Let $c \colon M \to S^n$ the collapse map which is a map of degree
one. Since $M$ is by assumption a homology sphere, it induces an
isomorphism on integral homology. By the Atiyah Hirzebruch
spectral sequence it induces isomorphisms
$$H_p(c;\bfL\langle 1 \rangle) \colon H_p(M;\bfL\langle 1 \rangle) \to H_p(S^n;\bfL\langle 1 \rangle)$$
for all $p \in \bbZ$. We obtain the following commutative diagram whose vertical arrows are parts of the
long exact surgery sequence~\eqref{Ranickis algebraic surgery sequence} where we
here use the decoration $s$, i.e. we take the Whitehead torsion into account.

$$\begin{CD}
H_{n+1}(M;\bfL\langle 1 \rangle) @> H_{n+1}(c;\bfL\langle 1 \rangle) >\cong > H_{n+1}(S^n;\bfL\langle 1 \rangle)
\\
@VVV @VVV
\\
L_{n+1}^s(\bbZ \pi) @> L_{n+1}^s(\pi_1(c)) >> L_{n+1}^s(\bbZ)
\\
@VVV @VVV
\\
\cals^{\topo}(M) @>\cals^{\topo}(c) >> \cals^{\topo}(S^n)
\\
@VVV @VVV
\\
H_n(M;\bfL\langle 1 \rangle) @> H_n(c;\bfL\langle 1 \rangle) > \cong > H_n(S^n;\bfL\langle 1 \rangle)
\\
@VVV @VVV
\\
L_n^s(\bbZ \pi) @> L_n^s(\pi_1(c)) >> L_n^s(\bbZ)
\end{CD}
$$
By the Poincar\'e Conjecture $\cals^{\topo}(S^n)$ is trivial.
An easy diagram chase together with Theorem~\ref{the:surgery_criterion_for_Borel}
shows that $M$ is strongly Borel if and only if
$$H_{n+1}(c;\bfL\langle 1 \rangle) \colon H_{n+1}(M;\bfL\langle 1 \rangle) \to L_{n+1}^s(\bbZ \pi)$$
is surjective. Notice that the natural map $\bfL^s \to \bfL$ is a homotopy equivalence since the
Whitehead group of the trivial group is zero. The following  diagram
$$\begin{CD}
H_{n+1}(\pt;\bfL^s\langle 1 \rangle) @> H_{n+1}(i;\bfL^s\langle 1 \rangle) >\cong>
H_{n+1}(M;\bfL^s\langle 1 \rangle) @> H_{n+1}(c;\bfL^s\langle 1 \rangle) >\cong >
H_{n+1}(S^n;\bfL^s\langle 1 \rangle).
\\
@VV \cong V @VVV @VV\cong V
\\
L_{n+1}^s(\bbZ) @> L_{n+1}^s(j) >>
L_{n+1}^s(\bbZ \pi) @> L_n^s(\pi_1(c)) >> L_{n+1}^s(\bbZ)
\end{CD}$$
commutes, where $i \colon \pt \to M$ and $j \colon \bbZ \to \bbZ
\pi$ denote the obvious inclusions and all maps marked with
$\cong$ are isomorphisms by the Atiyah-Hirzebruch spectral
sequence. Obviously $L_n^s(\pi_1(c)) \circ L_{n+1}^s(j)$ is the
identity. Hence $L_{n+1}^s(j)$ is bijective if and only if
$H_{n+1}(c;\bfL\langle 1 \rangle)$ is surjective. This shows that
$M$ is strongly Borel if and only if $L_{n+1}^s(j)$ is bijective.
\\[1mm]
\eqref{the:Integral_homology_spheres:nec_n+1}
Suppose that $M$ is strongly Borel and a rational homology sphere.
Since then $\cals^{\topo}(M) = \{0\}$, the map
$H_{n+1}(M;\bfL^s\langle 1 \rangle) \to L_{n+1}^s(\bbZ \pi)$ is surjective. Since this map
factorizes as
$$H_{n+1}(M;\bfL^s\langle 1 \rangle) \xrightarrow{H_{n+1}(f;\bfL^s\langle 1 \rangle)} H_{n+1}(B\pi;\bfL^s\langle 1 \rangle)
\to  L_{n+1}^s(\bbZ \pi)$$
for the classifying map $f \colon M \to B\pi$ and the latter map is rationally injective by assumption,
the homomorphism $H_{n+1}(f;\bfL^s\langle 1 \rangle)$ is rationally surjective.
Given a $CW$-complex $X$ we have the isomorphism given by the Chern character
$$\bigoplus_{i \ge 1} H_{n+1-4i}(M;\bbQ) \xrightarrow{\cong} H_{n+1}(M;\bfL^s\langle 1 \rangle).$$
Hence the map $H_{n+1-4i}(f;\bbQ) \colon H_{n+1-4i}(M;\bbQ)  \to H_{n+1-4i}(B\pi;\bbQ)$ is surjective for $i~\ge~1$.
\end{proof}

\typeout{-------------------- References -------------------------------}

\bibliographystyle{abbrv}
\bibliography{dbdef,dbpub,dbpre,dbklextra}

\def\cprime{$'$} \def\polhk#1{\setbox0=\hbox{#1}{\ooalign{\hidewidth
  \lower1.5ex\hbox{`}\hidewidth\crcr\unhbox0}}}
\begin{thebibliography}{10}

\bibitem{Adams(1966)}
J.~F. Adams.
\newblock On the groups ${J}({X})$. {I}{V}.
\newblock {\em Topology}, 5:21--71, 1966.

\bibitem{Bak(1981)}
A.~Bak.
\newblock {\em ${K}$-theory of forms}.
\newblock Princeton University Press, Princeton, N.J., 1981.

\bibitem{Browder(1969)}
W.~Browder.
\newblock The {K}ervaire invariant of framed manifolds and its generalization.
\newblock {\em Ann. of Math. (2)}, 90:157--186, 1969.

\bibitem{Cappell(1976a)}
S.~E. Cappell.
\newblock A splitting theorem for manifolds.
\newblock {\em Invent. Math.}, 33(2):69--170, 1976.

\bibitem{Chang-Weinberger(2003b)}
S.~Chang and S.~Weinberger.
\newblock On invariants of {H}irzebruch and {C}heeger-{G}romov.
\newblock {\em Geom. Topol.}, 7:311--319 (electronic), 2003.

\bibitem{Cohen(1973)}
M.~M. Cohen.
\newblock {\em A course in simple-homotopy theory}.
\newblock Springer-Verlag, New York, 1973.
\newblock Graduate Texts in Mathematics, Vol. 10.

\bibitem{Dold(1962)}
A.~Dold.
\newblock Relations between ordinary and extraordinary homology.
\newblock Colloq. alg. topology, Aarhus 1962, 2-9, 1962.

\bibitem{Farrell(1996)}
F.~T. Farrell.
\newblock {\em Lectures on surgical methods in rigidity}.
\newblock Published for the Tata Institute of Fundamental Research, Bombay,
  1996.

\bibitem{Farrell(2002)}
F.~T. Farrell.
\newblock The {B}orel conjecture.
\newblock In F.~T. Farrell, L.~G\"ottsche, and W.~L{\"u}ck, editors, {\em High
  dimensional manifold theory}, number~9 in ICTP Lecture Notes, pages 225--298.
  Abdus Salam International Centre for Theoretical Physics, Trieste, 2002.
\newblock Proceedings of the summer school ``High dimensional manifold theory''
  in Trieste May/June 2001, Number~1.
  http://www.ictp.trieste.it/\~{}pub\_off/lectures/vol9.html.

\bibitem{Farrell-Jones(1989)}
F.~T. Farrell and L.~E. Jones.
\newblock A topological analogue of {M}ostow's rigidity theorem.
\newblock {\em J. Amer. Math. Soc.}, 2(2):257--370, 1989.

\bibitem{Farrell-Jones(1990)}
F.~T. Farrell and L.~E. Jones.
\newblock Rigidity and other topological aspects of compact nonpositively
  curved manifolds.
\newblock {\em Bull. Amer. Math. Soc. (N.S.)}, 22(1):59--64, 1990.

\bibitem{Farrell-Jones(1993a)}
F.~T. Farrell and L.~E. Jones.
\newblock Isomorphism conjectures in algebraic ${K}$-theory.
\newblock {\em J. Amer. Math. Soc.}, 6(2):249--297, 1993.

\bibitem{Farrell-Jones(1993c)}
F.~T. Farrell and L.~E. Jones.
\newblock Topological rigidity for compact non-positively curved manifolds.
\newblock In {\em Differential geometry: Riemannian geometry (Los Angeles, CA,
  1990)}, pages 229--274. Amer. Math. Soc., Providence, RI, 1993.

\bibitem{Farrell-Jones(1998)}
F.~T. Farrell and L.~E. Jones.
\newblock Rigidity for aspherical manifolds with $\pi\sb 1\subset {GL}\sb
  m(\mathbb{{R}})$.
\newblock {\em Asian J. Math.}, 2(2):215--262, 1998.

\bibitem{Ferry-Ranicki-Rosenberg(1995)}
S.~C. Ferry, A.~A. Ranicki, and J.~Rosenberg.
\newblock A history and survey of the {N}ovikov conjecture.
\newblock In {\em Novikov conjectures, index theorems and rigidity, Vol.\ 1
  (Oberwolfach, 1993)}, pages 7--66. Cambridge Univ. Press, Cambridge, 1995.

\bibitem{Freedman(1983)}
M.~H. Freedman.
\newblock The disk theorem for four-dimensional manifolds.
\newblock In {\em Proceedings of the International Congress of Mathematicians,
  Vol.\ 1, 2 (Warsaw, 1983)}, pages 647--663, Warsaw, 1984. PWN.

\bibitem{Freedman-Quinn(1990)}
M.~H. Freedman and F.~Quinn.
\newblock {\em Topology of 4-manifolds}.
\newblock Princeton University Press, Princeton, NJ, 1990.

\bibitem{Hambleton-Kreck(1993c)}
I.~Hambleton and M.~Kreck.
\newblock Cancellation, elliptic surfaces and the topology of certain
  four-manifolds.
\newblock {\em J. Reine Angew. Math.}, 444:79--100, 1993.

\bibitem{Hempel(1976)}
J.~Hempel.
\newblock {\em $3$-{M}anifolds}.
\newblock Princeton University Press, Princeton, N. J., 1976.
\newblock Ann. of Math. Studies, No. 86.

\bibitem{Kirby-Siebenmann(1977)}
R.~C. Kirby and L.~C. Siebenmann.
\newblock {\em Foundational essays on topological manifolds, smoothings, and
  triangulations}.
\newblock Princeton University Press, Princeton, N.J., 1977.
\newblock With notes by J.~Milnor and M.~F.~Atiyah, Annals of Mathematics
  Studies, No. 88.

\bibitem{Kreck(1999)}
M.~Kreck.
\newblock Surgery and duality.
\newblock {\em Ann. of Math. (2)}, 149(3):707--754, 1999.

\bibitem{Kreck(2004b)}
M.~Kreck.
\newblock Cancellation for stable diffeomorphisms.
\newblock preprint, 2004.

\bibitem{Kreck-Lueck(2005)}
M.~Kreck and W.~L\"uck.
\newblock {\em The {N}ovikov Conjecture: Geometry and Algebra}, volume~33 of
  {\em Oberwolfach Seminars}.
\newblock Birkh\"auser, 2005.

\bibitem{Kreck-Leichtnam-Lueck(2002)}
E.~Leichtnam, W.~L{\"u}ck, and M.~Kreck.
\newblock On the cut-and-paste property of higher signatures of a closed
  oriented manifold.
\newblock {\em Topology}, 41(4):725--744, 2002.

\bibitem{Lueck(2002c)}
W.~L{\"u}ck.
\newblock A basic introduction to surgery theory.
\newblock In F.~T. Farrell, L.~G\"ottsche, and W.~L{\"u}ck, editors, {\em High
  dimensional manifold theory}, number~9 in ICTP Lecture Notes, pages 1--224.
  Abdus Salam International Centre for Theoretical Physics, Trieste, 2002.
\newblock Proceedings of the summer school ``High dimensional manifold theory''
  in Trieste May/June 2001, Number~1.
  http://www.ictp.trieste.it/\~{}pub\_off/lectures/vol9.html.

\bibitem{Lueck-Reich(2005)}
W.~L{\"u}ck and H.~Reich.
\newblock The {B}aum-{C}onnes and the {F}arrell-{J}ones {C}onjectures in {K}-
  and {L}-{T}heory.
\newblock In {\em {H}andbook of {$K$}-theory}, volume~2, pages 703 -- 842.
  Springer-Verlag, Berlin, 2005.

\bibitem{Milnor(1958b)}
J.~Milnor.
\newblock On simply connected {$4$}-manifolds.
\newblock In {\em Symposium internacional de topolog\'\i a algebraica
  International symposi um on algebraic topology}, pages 122--128. Universidad
  Nacional Aut\'onoma de M\'exico and UNESCO, Mexico City, 1958.

\bibitem{Ranicki(2005z)}
A.~Ranicki.
\newblock On the manifold structure set.
\newblock in preparation, 2005.

\bibitem{Ranicki(1973a)}
A.~A. Ranicki.
\newblock Algebraic {$L$}-theory. {I}. {F}oundations.
\newblock {\em Proc. London Math. Soc. (3)}, 27:101--125, 1973.

\bibitem{Ranicki(1980a)}
A.~A. Ranicki.
\newblock The algebraic theory of surgery. {I}. {F}oundations.
\newblock {\em Proc. London Math. Soc. (3)}, 40(1):87--192, 1980.

\bibitem{Ranicki(1981)}
A.~A. Ranicki.
\newblock {\em Exact sequences in the algebraic theory of surgery}.
\newblock Princeton University Press, Princeton, N.J., 1981.

\bibitem{Ranicki(1992)}
A.~A. Ranicki.
\newblock {\em Algebraic ${L}$-theory and topological manifolds}.
\newblock Cambridge University Press, Cambridge, 1992.

\bibitem{Serre(1953)}
J.-P. Serre.
\newblock Groupes d'homotopie et classes de groupes ab\'eliens.
\newblock {\em Ann. of Math. (2)}, 58:258--294, 1953.

\bibitem{Turaev(1988)}
V.~G. Turaev.
\newblock Homeomorphisms of geometric three-dimensional manifolds.
\newblock {\em Mat. Zametki}, 43(4):533--542, 575, 1988.
\newblock translation in Math. Notes 43 (1988), no. 3-4, 307--312.

\bibitem{Waldhausen(1978b)}
F.~Waldhausen.
\newblock Algebraic ${K}$-theory of topological spaces. {I}.
\newblock In {\em Algebraic and geometric topology (Proc. Sympos. Pure Math.,
  Stanford Univ., Stanford, Calif., 1976), Part 1}, pages 35--60. Amer. Math.
  Soc., Providence, R.I., 1978.

\bibitem{Wall(1962)}
C.~T.~C. Wall.
\newblock Classification of {$(n-1)$}-connected {$2n$}-manifolds.
\newblock {\em Ann. of Math. (2)}, 75:163--189, 1962.

\bibitem{Wall(1999)}
C.~T.~C. Wall.
\newblock {\em Surgery on compact manifolds}.
\newblock American Mathematical Society, Providence, RI, second edition, 1999.
\newblock Edited and with a foreword by A. A. Ranicki.

\end{thebibliography}

\noindent
Matthias Kreck\\
Mathematisches Institut, Universit\"at Heidelberg\\
Im Neuenheimer Feld 288, 69120 Heidelberg, Germany\\
kreck@mathi.uni-heidelberg.de\\
www.mathi.uni-heidelberg.de/$\sim$kreck/\\
Fax: 06221-545618\\[3mm]
Wolfgang L\"uck
\\Fachbereich Mathematik, Universit\"at M\"unster\\ 
Einsteinstr.~62, 48149 M\"unster, Germany\\
lueck@math.uni-muenster.de\\
http://www.math.uni-muenster.de/u/lueck\\
FAX: 49 251 8338370\\[3mm]

\end{document}